\documentclass[10pt,a4paper]{amsart}
\usepackage{graphicx,color}
\usepackage{url}

\usepackage{mathrsfs}

\newtheorem{theorem}{Theorem}[section]
\newtheorem{lemma}[theorem]{Lemma}           
\newtheorem{cor}[theorem]{Corollary}
\newtheorem{prop}[theorem]{Proposition}

\theoremstyle{definition}
\newtheorem{definition}[theorem]{Definition}

\theoremstyle{remark}

\numberwithin{equation}{section}

\makeatletter
\@namedef{subjclassname@2020}{%
  \textup{2020} Mathematics Subject Classification}
\makeatother

\subjclass[2020]{Primary 81U24, Secondary 47A75}
\keywords{quantum walk, eigenvalue, resonance}

\title[Complex translation methods and resonances for QWs]
{Complex translation methods and its application to resonances for quantum walks}
\author[K. Higuchi]{Kenta Higuchi}
\address[K. Higuchi]{JSPS Research Fellow PD, Graduate School of Science and Engineering, Ehime University, Bunkyo-cho 3, Matsuyama, Ehime, 790-8577, Japan}
\email{higuchi.kenta.en@ehime-u.ac.jp}
\author[H. Morioka]{Hisashi Morioka}
\address[H. Morioka]{Graduate School of Science and Engineering,
Ehime University, Bunkyo-cho 3, Matsuyama, Ehime, 790-8577, Japan}
\email{morioka@cs.ehime-u.ac.jp}

\thanks{K. Higuchi is supported by the JSPS Grant-in-aid for JSPS Fellow No. JP22J00430 and No. JP22KJ2364.
H. Morioka is supported by the JSPS Grant-in-aid for young scientists No. JP20K14327. }

\date{\today}

\begin{document}
\maketitle

\begin{abstract} 
In this paper, some properties of resonances for multi-dimensional quantum walks are studied.
Resonances for quantum walks are defined as eigenvalues of complex translated time evolution operators in the pseudo momentum space.
For some typical cases, we show some results of existence or nonexistence of resonances.
One is a perturbation of an elastic scattering of a quantum walk which is an analogue of classical mechanics.
Another one is a shape resonance model which is a perturbation of a quantum walk with a non-penetrable barrier.
\end{abstract}

%
%
\section{Introduction}
\begin{figure}[t]
\centering
\includegraphics[bb=0 0 808 300, width=10cm]{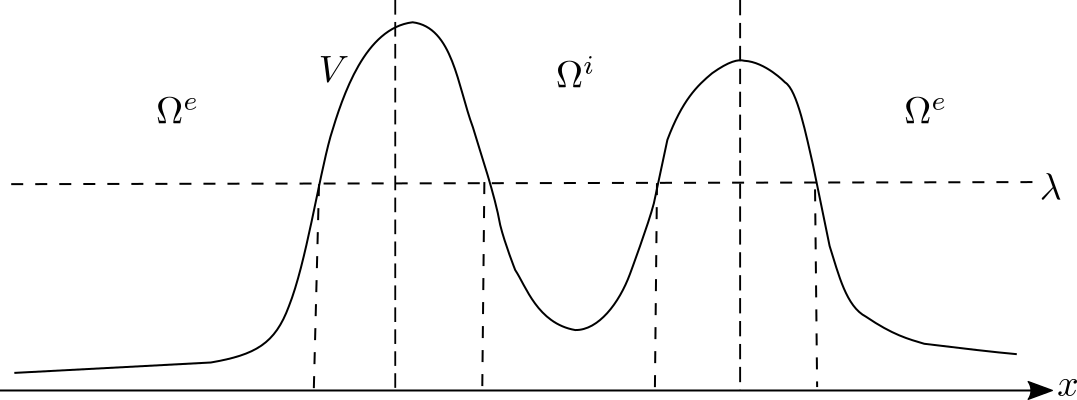}
\caption{The potential $V\in C^{\infty} ({\bf R}^d )$ for the shape resonance model $H(h)= -(h^2 /2)\Delta +V$ on ${\bf R}^d$ with small parameter $h>0$. If $V(x)\to 0$ rapidly as $|x| \to \infty$, $H(h)$ has no positive eigenvalue. However, there exist some resonances of $H(h)$ near the positive semi-axis.}
\label{fig_shaperesonance}
\end{figure}

In the research area of quantum physics, the study of resonances has a long history.
For Schr\"{o}dinger operators, one often adopts the framework of semi-classical analysis.
The shape resonance models (\cite{Kl}, \cite{CDKS}, \cite{Na1}, \cite{Na2}, \cite{Ka}) deserve our attention in this paper.
A typical shape resonance model is the Schr\"{o}dinger operator $H(h)=- (h^2 /2)\Delta +V$ on ${\bf R}^d$ with a small parameter $h>0$ where the potential $V$ is like a cut-off harmonic oscillator (see Figure \ref{fig_shaperesonance}).
In the classical mechanics, if a particle has an energy $\lambda$ in Figure \ref{fig_shaperesonance}, we can take a hyper-surface $K$ such that $K$ splits ${\bf R}^d$ into a connected bounded domain $\Omega^i$ and a connected exterior domain $\Omega^e$, and $K$ can be taken in the classically forbidden region $J(\lambda)=\{ x\in {\bf R}^d \ ; \ V(x)>\lambda\} $. 
Thus there are some trajectories bounded in the domain $\Omega^i$.
On the other hand, the Schr\"{o}dinger operator has no bound state since there is the tunneling effect in quantum mechanics.
Combes et al. \cite{CDKS} showed that there exist some resonances of $ H(h) $ in a small neighborhood of Dirichlet eigenvalues of $-(h^2 /2)\Delta +V$ in $\Omega^i$ if we take sufficiently small $h>0$.
Klein \cite{Kl} proved the absence of resonances in a neighborhood of the essential spectrum of $H(h)$ (which is the positive-semi axis) for the case where $V$ is a non-trapping potential.
Their results suggest that the existence of resonances near the positive semi-axis reflects the existence of bounded classical trajectories as well as scattering trajectories. 
We would like to refer \cite{Si} and \cite{DyZw} for general information.
See also the references therein.

In these previous works, the complex dilation (\cite{AgCo}) for Schr\"{o}dinger operators was often used.
The essential spectrum is rotated and becomes a half line in the lower half plane.
On the other hand, eigenvalues and resonances are invariant under the complex dilation after the deformed essential spectrum moves over them.
Thus the complex dilation allows us to study resonances as isolated eigenvalues of the dilated Hamiltonian.

In this paper, we study resonances in the context of quantum walks (QWs) as a perturbation problem of eigenvalues of time evolution operators.
The spectral theory and the eigenvalue problem for time evolution operators of QWs on ${\bf Z}^d$ are studied in some recent works.
The properties of the essential spectrum are studied in \cite{RST}, \cite{Mo}, \cite{MoSe}, \cite{Ti}, \cite{KKMS}, \cite{KKMS2}.
As an earlier work than these articles, Kato-Kuroda \cite{KaKu} presented an abstract theory of wave operators in view of perturbations of unitary operators rather than that of self-adjoint operators.
Their theory corresponds to considering the discrete unitary group associated with the unitary operators rather than the continuous one-parameter unitary group associated with self-adjoint operators.

We introduce the model of the $d$-dimensional QW which is studied in this paper.
For the sake of simplicity of notations, we restrict our consideration to the case $d=2$.
For other $d$, the argument is parallel.
In this paper, we consider finite rank perturbations of the free QW.
Let $ e_1 = [1,0]^{\mathsf{T}}$ and $e_2 = [0,1]^{\mathsf{T}} $.
The free QW is given by the unitary operator $U_0 =S$ on $\mathcal{H}:= \ell^2 ({\bf Z}^2 ; {\bf C}^4 )$ where $S$ is the shift operator
$$
(Su)(x)= \begin{bmatrix} 
u _{\leftarrow} (x+e_1 ) \\ u_{\rightarrow} (x-e_1 ) \\ u_{\downarrow} (x+e_2 ) \\ u_{\uparrow} (x-e_2 ) \end{bmatrix} , \quad x\in {\bf Z}^2 ,
$$
for a ${\bf C}^4 $-valued sequence $ u=\{ [u_{\leftarrow} (x), u_{\rightarrow} (x), u_{\downarrow} (x), u_{\uparrow} (x) ]^{\mathsf{T}} \} _{x\in {\bf Z}^2} $.
The perturbed QW is defined by a unitary operator $U=SC$ on $\mathcal{H}$ where $C$ is the operator of multiplication by a matrix $ C(x)\in \mathrm{U} (4)$ at every $x\in {\bf Z}^2$.
Note that the time evolutions of the free QW and the perturbed QW are define by
$$
\Psi (t,\cdot )= U^t \psi , \quad \Psi_0 (t,\cdot )= U_0^t \psi , \quad t\in {\bf Z} ,
$$
for an initial state $\psi $.
If $\psi \in \mathcal{H}$, we have 
$$
\| \Psi (t,\cdot ) \| _{\mathcal{H}} = \| \Psi_0 (t,\cdot ) \| _{\mathcal{H}} = \| \psi \| _{\mathcal{H}} ,
$$
for any $t\in {\bf Z} $ since $U$ and $U_0$ are unitary on $\mathcal{H}$.
Throughout this paper, we assume the following property.
Due to this assumption, the perturbation $V:= U-U_0$ is of finite rank.

\medskip

{\bf (A-1)} There exists a positive integer $ M_0 $ such that $ C(x)=I_4$ which is the $4\times 4$ identity matrix for $x\in {\bf Z}^2 \setminus \Omega^i$ where $\Omega^i = \{ x\in {\bf Z}^2 \ ; \ |x_1| \leq M_0, \ |x_2| \leq M_0 \} $.

\medskip

In general, the operator $U$ may have some eigenvalues. 
In the settings introduced as above, associated eigenstates are localized in the domain $\Omega^i$.
Under some suitable conditions, we can find another QW $U'$ such that $U'$ has no eigenvalue even if $U'-U$ is sufficiently small in a suitable topology and $U$ has some eigenvalues.
This kind of situations motivates us to study resonances of QWs.
Namely, the eigenvalues of $U$ may move into the second sheet of the Riemann surface.
In order to define resonances of the QW $U$, we consider the meromorphic extension of the resolvent operator $R(\kappa)=(U-e^{-i\kappa} )^{-1} $ with $\mathrm{Im} \, \kappa >0$ to the lower half region of the complex torus ${\bf T} _{{\bf C}}: = {\bf C} /2\pi {\bf Z}$ as the second sheet of the Riemann surface.
In this paper, we use the complex translation in the pseudo momentum space as an analogue of the complex dilation.
We note that a resonance expansion in view of the dynamics of one-dimensional QWs is going to be shown in the forthcoming paper \cite{HMS}.

This paper is organized as follows.
At the beginning of Section 2, we recall some known facts in the spectral theory for QWs.
After that, we introduce the complex translation in the pseudo momentum space.
The rigorous definition of (outgoing) resonances is given here.
In Section 3, we study elastic scattering of QWs.
Some relations between eigenvalues or resonances and closed trajectories are given here. 
Elastic scattering of QWs is an analogue of classical mechanics.
After that, we also study a QW with a non-penetrable barrier on the boundary of $\Omega^i$.
This QW is split into a QW on a finite graph and on its exterior domain.
Since the exterior QW is a finite perturbation of the free QW, its spectrum consisits of the essential spectrum on the unit circle.
Thus some eigenvalues which are embedded in the essential spectrum appear in this model due to the QW on the finite graph.
We regard this model as an analogue of the Schr\"{o}dinger operator with the Dirichlet boundary condition, since the spectrum of the interior Dirichlet Schr\"{o}dinger operator consists of discrete eigenvalues.
We prepare some estimates of resolvent operators.
In Section 4, an example of resonances is given in a constructive approach.
Proposition \ref{S4_prop_ecperturbation} is our result of this section.
Namely, we consider $U$ for the case where $U$ has (approximately) a few closed trajectories in view of the elastic scattering.
Next we consider an analogue of shape resonance models in Section 5.
Here the non-penetrable barrier is perturbed and the perturbed barrier is the simplest case of analogues for the classically forbidden region $J(\lambda)$ in the standard quantum mechanics. 
We prove that some eigenvalues of a QW with non-penetrable barrier move into the second sheet of the Riemann surface due to the perturbation.
As a conclusion, we show Theorem \ref{S5_thm_existenceres}.
Corollary \ref{S5_cor_main} is our main result.

The notations which are used throughout this paper are as follows.
As above, we put $ \mathcal{H} = \ell^2 ({\bf Z}^2 ; {\bf C}^4 )$ equipped with the inner product
$$
(f,g) _{\mathcal{H}} = \sum _{x\in {\bf Z}^2} \sum _{j\in \{ \leftarrow , \rightarrow , \downarrow , \uparrow \} } f_j (x) \overline{g_j (x)} , \quad f,g\in \mathcal{H}.
$$
We often use the standard basis on the vector space ${\bf R}^2 $ and ${\bf C}^4$.
We denote it by 
$$
e_1 = [1,0]^{\mathsf{T}} , \quad e_2 = [0,1]^{\mathsf{T}} ,
$$
 and 
 $${\bf e} _{\leftarrow} = [1,0,0,0]^{\mathsf{T}} , \quad {\bf e} _{\rightarrow} = [0,1,0,0]^{\mathsf{T}} ,\quad {\bf e} _{\downarrow} = [0,0,1,0] ^{\mathsf{T}} ,\quad {\bf e} _{\uparrow} = [0,0,0,1]^{\mathsf{T}} .
 $$
For an operator $A$ on $ \mathcal{H} $, we denote by $\sigma (A)$, $\sigma_{ess} (A)$, $\sigma_p (A)$ and $\sigma_{ac} (A)$ the spectrum, the essential spectrum, the point spectrum and the absolutely continuous spectrum of $A$, respectively.
For Banach spaces $ \mathcal{H}_1 $ and $ \mathcal{H}_2 $, ${\bf B} ( \mathcal{H}_1 ; \mathcal{H}_2 ) $ denotes the space of bounded linear operators from $ \mathcal{H}_1 $ to $ \mathcal{H}_2 $.
If $ \mathcal{H}_1 = \mathcal{H}_2 $, we simply write ${\bf B} (\mathcal{H}_1 ; \mathcal{H}_1 )= {\bf B} (\mathcal{H}_1 )$.
The flat torus is defined by ${\bf T}^d = {\bf R}^d / 2\pi {\bf Z}^d $.
The complex torus is defined by ${\bf T} _{{\bf C}}^d = {\bf C}^d /2\pi {\bf Z} ^d = {\bf T} ^d +i{\bf R} ^d $.
For $a\in {\bf R}$, we put 
$$
\mathcal{O}^+_a = \{ \kappa \in {\bf T} _{{\bf C}} \ ; \ \mathrm{Im} \, \kappa >a \} , \quad \mathcal{O}^-_a = \{ \kappa \in {\bf T} _{{\bf C}} \ ; \ \mathrm{Im} \, \kappa < a \} .
$$

\section{Complex translation on the pseudo momentum space}
\subsection{Preliminary results of spectra}
First of all, let us recall the spectra of $U_0$ and $U$ under the assumption (A-1).
The Fourier transform is defined by 
$$
\widehat{u} (\xi )=(\mathcal{F}u)(\xi )= \frac{1}{2\pi} \sum _{x\in {\bf Z}^2} e^{-ix\cdot \xi} u(x), \quad \xi \in {\bf T}^2 ,
$$
where $x\cdot \xi = x_1 \xi_1 +x_2 \xi_2 $.
It is well-known that $\mathcal{F}$ is a unitary operator from $\mathcal{H}$ to $\widehat{\mathcal{H}} := L^2 ({\bf T}^2 ; {\bf C}^4 ) $.
The operator $ \widehat{U}_0 = \mathcal{F} U_0 \mathcal{F}^* $ is multiplication by the $4\times 4$ unitary matrix
$$
\widehat{U}_0 (\xi )= \mathrm{diag} [ e^{i\xi_1} , e^{-i\xi_1} , e^{i\xi_2} , e^{-i\xi_2} ] .
$$
For $\kappa \in {\bf C}$, we have 
$$
\det (\widehat{U}_0 (\xi ) - e^{-i\kappa} ) = (e^{i\xi_1} -e^{-i\kappa} )(e^{-i\xi_1} -e^{-i\kappa} )(e^{i\xi_2} -e^{-i\kappa} )(e^{-i\xi_2} -e^{-i\kappa} ).
$$
This formula determines the spectrum of $U_0$.

\begin{lemma}
We have $\sigma (U_0)= \sigma_{ac} (U_0)= \{ e^{-i\lambda} \ ; \ \lambda \in [0,2\pi )\} =S^1 $. 
\label{S2_lem_specU0}
\end{lemma}

Since $ V=U-U_0$ is of finite rank, we can see that the essential spectrum of $U$ coincides with $\sigma_{ess} (U_0)=S^1$.
For details of this topic, the rigorous proof was given in \cite{MoSe} which is parallel to the well-known Weyl's singular sequence lemma for compact perturbations of self-adjoint operators (see e.g. \cite{Ya}).
We can also prove the absence of the singular continuous spectrum of $U$ (see \cite[Theorem 2.4]{RST} for 1DQWs and \cite[Lemma 4.16]{KKMS2} for 2DQWs).
If there exist some eigenvalues, they are embedded in the essential spectrum.

\begin{lemma}
Under the assumption (A-1), we have $ \sigma_{ess} (U)=S^1$ and $\sigma_p (U)$ consists of eigenvalues lying on $S^1$ with finite multiplicities.
\label{S2_lem_specU}
\end{lemma}

\textit{Remark.}
If we add another condition (C), we can show the absence of eigenvalues (\cite{MoSe}, \cite{KKMS}, \cite{MS4} for 1DQWs, and \cite{KKMS2} for 2DQWs).

\medskip

{\bf (C)} Let $C(x)=[ c_{j,k} (x) ] _{j,k\in \{ \leftarrow , \rightarrow , \downarrow , \uparrow \} }$ for every $x\in \Omega^i$.
One of the following properties holds true.
\begin{enumerate}
\item
Neither $\det [c_{j,k} (x)] _{j,k\in \{ \leftarrow , \downarrow \} }$ nor $\det [c_{j,k} (x)] _{j,k\in \{ \rightarrow , \uparrow \} } $ vanishes for any $x\in \Omega^i$.

\item
Neither $\det [c_{j,k} (x)] _{j,k\in \{ \leftarrow , \uparrow \} }$ nor $\det [c_{j,k} (x)] _{j,k\in \{ \rightarrow , \downarrow \} } $ vanishes for any $x\in \Omega^i$.
\end{enumerate}

\medskip

Note that the condition (C) is sufficient for the absence of eigenvalues.
For the most part of our argument, we do not assume that (C) holds true.
We use (C) in a typical case appearing in the context of scattering theory.

\subsection{Complex translation on the pseudo momentum space}
Suppose that an operator $A$ on a Hilbert space $\mathbb{H}$ has an isolated eigenvalue $\lambda \in {\bf C}$.
The algebraic multiplicity of $\lambda$ is defined by the rank of the operator 
$$
\widetilde{P}_A (\lambda)=-\frac{1}{2\pi i} \oint _{\widetilde{\mathcal{L}}(\lambda)} (A-z)^{-1} dz,
$$
where $\widetilde{\mathcal{L}} (\lambda )$ is a sufficiently small counterclockwise loop without self-intersection such that there is no other eigenvalues inside $\widetilde{\mathcal{L}} (\lambda )$.
Then the resolvent operator $(A-z)^{-1}$ acts a crucial role in the study of eigenvalues.
On the other hand, if an eigenvalue $\lambda$ is embedded in the continuous spectrum of $A$, the projection $\widetilde{P}_A (\lambda)$ is not well-defined.
For the case $ \mathbb{H}=\mathcal{H}$ and $A=U$, eigenvalues are embedded in the essential spectrum in view of Lemma \ref{S2_lem_specU}.
Then we have to avoid the difficulty for $\widetilde{P}_A (\lambda)$.
In order to do this, we deform the continuous spectrum $\sigma_{ess} (U)$ by using a complex translation on the pseudo momentum space ${\bf T}^2 $.

Let us begin with the real translation operator.
For $ \theta \in {\bf T} $, we define the operator $\widehat{T} (\theta )$ of translation by 
\begin{equation}
(\widehat{T} (\theta ) \widehat{u} )(\xi )= \begin{bmatrix} \widehat{u} _{\leftarrow} (\xi_1 - \theta , \xi_2 ) \\ \widehat{u} _{\rightarrow} (\xi_1 + \theta , \xi_2 ) \\ \widehat{u} _{\downarrow} (\xi_1 , \xi_2 - \theta ) \\ \widehat{u}_{\uparrow} (\xi_1 , \xi_2 + \theta ) \end{bmatrix} , \quad \widehat{u} \in L^2 ({\bf T}^2 ; {\bf C}^4 ).
\label{S2_eq_realtranslation}
\end{equation}
Obviously, $\widehat{T} ( \theta ) $ is unitary on $L^2 ({\bf T}^2 ; {\bf C}^4 )$ for every $\theta \in {\bf T} $.
We immediately see 
\begin{gather*}
\begin{split}
&\widehat{U}_0 (\theta ):= \widehat{T} ( \theta ) \widehat{U}_0 \widehat{T} (\theta )^{-1} = e^{-i\theta} \widehat{U}_0 , \\ 
&\widehat{R}_0 ( \kappa , \theta ):= \widehat{T} ( \theta ) \widehat{R}_0 (\kappa) \widehat{T} (\theta )^{-1} = e^{i\theta} \widehat{R}_0 (\kappa -\theta) , \quad \kappa \in {\bf T} _{{\bf C}} ,
\end{split}
\end{gather*}
 where
$$
R_0 (\kappa )=(U_0 -e^{-i\kappa} )^{-1} , \quad \widehat{R}_0 (\kappa )=\mathcal{F} R_0 (\kappa) \mathcal{F}^* .
$$
Thus $\widehat{U}_0 (\theta)$ as well as $\widehat{R}_0 (\kappa ,\theta)$ are naturally defined for $\theta \in \mathcal{O}^-_0 $ and $\kappa \in \mathcal{O}_{\mathrm{Im} \, \theta}^+$, and they are operators belonging to ${\bf B} (L^2 ({\bf T}^2 ; {\bf C}^4))$.
Precisely, we have 
$$
\widehat{U}_0 (\theta )= e^{-i \mathrm{Re} \, \theta} e^{\mathrm{Im} \, \theta} \widehat{U}_0 , \quad \theta \in \mathcal{O}^-_0 .
$$
As a direct consequence, the spectrum of $\widehat{U}_0$ is deformed as follows.

\begin{lemma}
We have $\sigma (\widehat{U}_0 (\theta ))= \sigma_{ac} (\widehat{U}_0 (\theta)) = \{ e^{\mathrm{Im} \, \theta} e^{i\lambda} \ ; \ \lambda \in [0,2\pi ) \} $ for $\theta \in \mathcal{O}^-_0$.
\label{S2_lem_specU0theta}
\end{lemma}

Let us turn to the complex translation $U(\theta)$ of $U$.
At first we consider the real translation again.
In order to construct $U(\theta)$, it is convenient to use $ T(\theta)= \mathcal{F}^* \widehat{T} (\theta ) \mathcal{F}$ for $\theta \in {\bf T}$.
In fact, we have 
\begin{equation}
(T(\theta)u)(x)= \mathrm{diag} [e^{i\theta x_1} , e^{-i\theta x_1} , e^{i \theta x_2 } , e^{-i \theta x_2} ] u(x) , \quad x\in {\bf Z}^2 ,
\label{S2_eq_complextranslation_lattice}
\end{equation}
for $\theta \in {\bf T}$.
We put
$$
V(\theta ):= T(\theta )VT(\theta )^{-1} .
$$
In view of the assumption (A-1), we can extend $V(\theta)$ to $\theta \in \mathcal{O}^-_0$ as a finite rank operator on $\mathcal{H}$.
Then we obtain a construction of the complex translation
$$
U(\theta)=T(\theta)UT(\theta)^{-1} =U_0 (\theta) +V(\theta ) \in {\bf B} (\mathcal{H}), \quad \theta \in \mathcal{O}^-_0,
$$
and formally define
$$
R(\kappa , \theta)= T(\theta)R(\kappa)T(\theta)^{-1} = (U(\theta)-e^{-i\kappa} )^{-1} \in {\bf B} (\mathcal{H} ), \quad e^{-i\kappa} \not\in \sigma (U(\theta)) .
$$
The compactness of $V(\theta )$ and the definition of $T(\theta)$ imply the following property of the continuous spectrum.

\begin{lemma}
For $\theta \in \mathcal{O}^-_0$, we have $ \sigma_{ess} (U(\theta)) = \sigma_{ess} (U_0 (\theta)) = \{ e^{\mathrm{Im} \, \theta} e^{i\lambda} \ ; \ \lambda \in [0,2\pi ) \} $.
\label{S2_lem_essspecUtheta}
\end{lemma}

Now let us consider the meromorphic extension of $R(\kappa)$ from $\mathcal{O}^+_0$ to the second sheet of the Riemann surface. We show that this meromorphic extension coincides with $R(\kappa,\theta)$ for $\theta \in \mathcal{O}^- _0 $ and $\kappa \in \mathcal{O}^+ _{\mathrm{Im} \, \theta} $. 
We put
$$
\widehat{\mathcal{D}} = \left\{ \widehat{f} | _{{\bf T}^2} \ ; \ \widehat{f} \text{ is a } {\bf C}^4 \text{-valued analytic function on } {\bf T}^2 _{{\bf C}} \right\} ,
$$
and
$$
\mathcal{D} = \mathcal{F}^* \widehat{\mathcal{D}} .
$$
Note that $\mathcal{D}$ is the set of ${\bf C}^4$-valued, super-exponentially decreasing sequences on ${\bf Z}^2$ in view of Paley-Wiener's theorem (\cite{Ves}) : 
$$
\mathcal{D} = \cap _{\gamma >0} \{ f\in \mathcal{H} \ ; \ e^{\tau \langle \cdot \rangle}f \in \mathcal{H} \text{ for any } 0\leq \tau <\gamma \} .
$$
The operator $T(\theta)$ for $\theta \in \mathcal{O}^-_0$ has $\mathcal{D}$ as its dense domain in $\mathcal{H} $.

\begin{lemma}
We define the function $F_f ( \kappa )$ for $\kappa \in \mathcal{O}_0^+$ and $f\in\mathcal{D}$ by
$$
F_f (\kappa )= (R(\kappa)f,f) _{\mathcal{H}} .
$$
For any $f\in \mathcal{D}$ and any $a<0$, the function $F_f (\kappa)$ has the meromorphic extension to $\mathcal{O}^+_{a} $ with poles of finite rank.
Poles of $F_f (\kappa) $ lie on ${\bf T}$ or in $\mathcal{O}^-_0 \cap \mathcal{O}^+ _a$.
\label{S2_lem_meromorphicextenstion}
\end{lemma}

Proof.
In this proof, we sometimes use the analytic Fredholm theory (see e.g. \cite[Theorem C.8]{DyZw}).

Since $U$ is a unitary operator on $\mathcal{H}$, $R(\kappa)$ is analytic in $\mathcal{O}^+_0 $ with respect to $\kappa$.
Thus $F_f (\kappa)$ with $\kappa \in \mathcal{O}^+_0 $ is well-defined for any $f\in \mathcal{H}$.
In view of the resolvent equation
\begin{equation}
R(\kappa)=R_0 (\kappa)-R(\kappa)VR_0 (\kappa)= R_0 (\kappa)-R_0 (\kappa)VR(\kappa),
\label{S2_eq_resolventinverse00}
\end{equation}
we have 
\begin{equation}
(1-VR(\kappa))(1+VR_0 (\kappa))=(1+VR_0 (\kappa))(1-VR(\kappa))=1.
\label{S2_eq_resolventinverse01}
\end{equation}
This implies that $1+VR_0 (\kappa) $ is invertible for $\kappa \in \mathcal{O}^+_0 $ and we obtain
\begin{equation}
R(\kappa)= R_0 (\kappa) (1+VR_0 (\kappa))^{-1} , \quad \kappa \in \mathcal{O}^+_0 .
\label{S2_eq_resolventinverse11}
\end{equation}
By the same way, we also see 
\begin{gather}
R(\kappa,\theta)=R_0 (\kappa,\theta)-R(\kappa,\theta) V (\theta )R_0 (\kappa,\theta)=R_0 (\kappa,\theta)-R_0 (\kappa,\theta) V (\theta )R (\kappa,\theta),
\label{S2_eq_resolventinverse02} \\ 
\begin{split}
&(1-V(\theta)R(\kappa,\theta))(1+V(\theta)R_0 (\kappa,\theta))
= (1+V(\theta)R_0 (\kappa,\theta))(1-V(\theta)R(\kappa,\theta))\\
&=1,
\end{split}
\label{S2_eq_resolventinverse03} \\
R(\kappa,\theta)= R_0 (\kappa,\theta) (1+V(\theta)R_0 (\kappa,\theta))^{-1} , 
\label{S2_eq_resolventinverse12}
\end{gather}
for $\kappa \in \mathcal{O}^+_0 $ and $\theta \in {\bf T}$.
Since $T(\theta )$ is unitary on $\mathcal{H}$ for $\theta \in {\bf T}$, we have 
\begin{gather}
\begin{split}
F_f (\kappa ) &= (T(\theta)R(\kappa )f,T(\theta)f)_{\mathcal{H}} \\
&= (R(\kappa,\theta) T(\theta)f,T(\theta)f)_{\mathcal{H}} \\
&= ( R_0 (\kappa ,\theta) (1+ V (\theta) R_0 (\kappa,\theta))^{-1} T( \theta )f, T(\theta )f)_{\mathcal{H}} ,
\end{split}
\label{S2_eq_mero00}
\end{gather}
for $\kappa \in \mathcal{O}^+_0 $.
Here we have used the formulas (\ref{S2_eq_resolventinverse00})-(\ref{S2_eq_resolventinverse12}).

We fix $ \kappa_0 \in \mathcal{O}^+ _0 $ and $f\in \mathcal{D}$.
If $ \theta \in \mathcal{O}^- _{\mathrm{Im} \, \kappa_0} $, we have $\mathrm{Im} \, (\kappa _0 -\theta)>0 $.
Then $1+ V (\theta) R_0 (\kappa_0,\theta) = e^{i\theta} V(\theta) R_0 (\kappa_0 - \theta)$ is compact and analytic with respect to $\theta \in \mathcal{O}^-_{\mathrm{Im} \, \kappa_0} $.
The existence of the inverse $(1+V(\theta)R_0 (\kappa_0 ,\theta ))^{-1}  $ for $\theta \in {\bf T} \subset \mathcal{O}^-_{\mathrm{Im} \, \kappa_0} $ follows from (\ref{S2_eq_resolventinverse03}).
The analytic Fredholm theory shows the existence of the meromorphic extension of $(1+ V (\theta) R_0 (\kappa_0,\theta) )^{-1} $ to $\theta \in \mathcal{O}^- _{\mathrm{Im} \, \kappa_0} $ with poles of finite rank.
Then the function 
$$
G_{f} (\kappa_0 ,\theta )=  ( R_0 (\kappa_0 , \theta) (1+ V (\theta) R_0 (\kappa_0,\theta) )^{-1} T(\theta) f, T(\theta)f )_{\mathcal{H} } ,
$$
of $\theta$ is meromorphic in $\mathcal{O}^-_{\mathrm{Im} \, \kappa_0} $.
Furthermore, the equality (\ref{S2_eq_mero00}) for $\theta \in {\bf T} $ implies that $G_{f} (\kappa_0 , \theta )=F_f (\kappa_0 )$ for any $\theta \in {\bf T} $.
Then the function $ G_{f} (\kappa_0 , \theta )$ is a constant in $\mathcal{O}^- _{\mathrm{Im} \, \kappa_0} $, i.e.,
$$
G_{f} (\kappa_0 , \theta )=F_f (\kappa_0) , \quad \theta \in \mathcal{O}^- _{\mathrm{Im} \, \kappa_0} .
$$
As a consequence, we see that $(1+ V (\theta) R_0 (\kappa_0,\theta) )^{-1} $ is analytic with respect to $\theta \in \mathcal{O}^- _{\mathrm{Im} \, \kappa_0} $.

Next we fix $\theta_0 \in \mathcal{O}^-_0 $ and $f\in \mathcal{D}$.
From the above argument, we have $ G_{f} (\kappa ,\theta_0 )= F_f (\kappa )$ for $\kappa \in \mathcal{O}^+_{0}$.
Noting $\mathrm{Im} \, \kappa - \mathrm{Im} \, \theta_0 >0 $ for $\kappa \in \mathcal{O}^+ _{\mathrm{Im} \, \theta_0} $, we see that $R_0 (\kappa , \theta_0 )= e^{i\theta_0} R_0 (\kappa - \theta_0 ) $ and $V(\theta_0 ) R_0 (\kappa, \theta_0 )=  e^{i\theta_0} V( \theta_0 ) R_0 ( \kappa - \theta_0 )$ are analytic with respect to $\kappa\in \mathcal{O}^+_{\mathrm{Im} \, \theta_0}$.
We also see that $V(\theta_0 ) R_0 (\kappa, \theta_0 )$ is compact for $\kappa\in \mathcal{O}^+_{\mathrm{Im} \, \theta_0}$.
The existence of the inverse $ (1+  V(\theta_0 ) R_0 (\kappa , \theta_0 ))^{-1} $ for some $\kappa\in \mathcal{O}^+_0 $ has been shown in the above argument, since $ (1+  V(\theta ) R_0 (\kappa , \theta ))^{-1} $ is analytic with respect to $\theta \in \mathcal{O}^- _{\mathrm{Im} \, \kappa} $.
Then the analytic Fredholm theory implies that $ (1+  V(\theta_0 ) R_0 (\kappa , \theta_0 ))^{-1}$ has the meromorphic extension to $\kappa \in \mathcal{O}^+ _{\mathrm{Im} \, \theta_0} $ with poles of finite rank.
Therefore, we obtain the meromorphic extension of $ F_f (\kappa)$ to the domain $\mathcal{O}^+ _{\mathrm{Im} \, \theta_0} $ with poles of finite rank.
Since $ \theta_0 \in \mathcal{O}^-_0 $ is arbitrary, we obtain the meromorphic extension of $F_f (\kappa )$ to $\mathcal{O}^+_{a} $ for any $a<0$.

Finally, we note that the discrete eigenvalues of $ U(\theta)$ in the region ${\bf C}^e _{\theta} := \{ z\in {\bf C} \ ; \ |z|> e^{\mathrm{Im} \, \theta} \} $ are invariant with respect to $\theta $.
An eigenvalue $ \lambda (\theta )$ of $U(\theta )$ moves analytically with respect to $\theta \in \mathcal{O}^+_{a} $ for any $a<0$ as long as $\lambda (\theta )$ belongs to ${\bf C}^e _{\theta} $.
For any $\alpha \in {\bf T}$, $U( \theta + \alpha ) $ and $ U(\theta )$ are unitary equivalent.
Then we have $\lambda ( \theta + \alpha )=\lambda (\theta )$.
By the analyticity of $\lambda (\theta )$, this implies that $\lambda (\theta )$ is a constant with respect to $\theta$ as long as $\lambda (\theta)\in {\bf C}^e _{\theta}$.
\qed

\medskip

In view of Lemma \ref{S2_lem_meromorphicextenstion} and its proof, the rigorous definition of resonances is given as follows.
Note that the resonances are independent of choice of $\theta$ as long as they are in ${\bf C}^e _{\theta} $ even though they are defined as eigenvalues of $ U( \theta )$.

\begin{definition}
Suppose that $\theta \in \mathcal{O}^-_0 $ is fixed.
We call $e^{-i\kappa} $ for $ \kappa \in \mathcal{O}^+ _{\mathrm{Im} \, \theta} $ a resonance of $U$ if $e^{-i\kappa} \in \sigma_p (U(\theta )) $.

\label{S2_def_resonance}
\end{definition}

Now we can define the algebraic multiplicity of a resonance $\lambda $ of $U$ by 
$$
\widetilde{P}_U ( \lambda )=- \frac{1}{2\pi i} \oint _{\widetilde{\mathcal{L}} (\lambda )} (U(\theta)-z )^{-1} dz,
$$
where $\widetilde{\mathcal{L}} (\lambda )$ is a sufficiently small counterclockwise loop without self-intersection and the inside of $\widetilde{\mathcal{L}} (\lambda )$ contains no eigenvalues except for $\lambda$.
However, it is convenient to replace $P_U (\kappa):= \widetilde{P}_U (e^{-i\kappa} )$ with $z=e^{-i\mu }$ by the integral of $\mu$ as 
\begin{equation}
P_U (\kappa )= \frac{1}{2\pi} \oint _{\mathcal{L} (\kappa)} e^{-i\mu} R(\mu ,\theta ) d\mu,
\label{S2_eq_algebraicmulti_proj}
\end{equation}
where $ \mathcal{L} (\kappa )$ is a sufficiently small counterclockwise loop without self-intersection such that there is no other poles inside $ \mathcal{L} (\kappa )$.
Note that we can take $\widetilde{\mathcal{L} } (e^{-i\kappa})$ such that $\mathcal{L} (\kappa)$ satisfies the above condition under the change of variable $z=e^{-i\mu}$.
See Lemma \ref{S3_lem_estimateRi}.

\begin{definition}
Suppose that $\theta \in \mathcal{O}^-_0 $ is fixed.
If $e^{-i\kappa} $ with $\kappa \in \mathcal{O}^+ _{\mathrm{Im} \, \theta} $ is a resonance of $U$, the algebraic multiplicity of $e^{-i\kappa}$ is defined by $ \mathrm{rank} \, P_U (\kappa ) $. 
\label{S2_def_multiresonance}
\end{definition}

\subsection{Outgoing solution}
If $e^{-i\lambda}$ with $\lambda \in {\bf T}$ is an eigenvalue of $U$, then $e^{-i\lambda}$ is an eigenvalue of $ U(\theta )$ for any $\theta \in \mathcal{O}^-_0 $ as we have seen in the proof of Lemma \ref{S2_lem_meromorphicextenstion}.
More precisely, the following property holds true.

\begin{lemma}
If $e^{-i\kappa} \in \sigma_p (U) $ with $\kappa \in {\bf T}$, then we have $e^{-i\kappa} \in \sigma_p (U(\theta )) $.
Each associated eigenfunction $u\in \mathcal{H}$ is supported in $ \Omega^i $.
\label{S2_lem_evUinvariant}
\end{lemma}

Proof.
It suffices to show $\mathrm{supp} \, u \subset \Omega^i $ for any eigenfunction $u$.
Suppose that $Uu=e^{-i\kappa} u$ holds true.
For $x\in {\bf Z}^2 \setminus \Omega^i$ with $x_1 \leq -M_0 -1$, we have $u_{\leftarrow} (x )=e^{-i\kappa} u_{\leftarrow} (x-e_1 )$.
Applying this equation at points $x,x-e_1,\ldots,x-(N-1)e_1 $ for any positive integer $N$, we have $ u_{\leftarrow} (x)= e^{-i\kappa N} u _{\leftarrow} (x-Ne_1 )$.
In view of $u\in \mathcal{H}$, $u_{\leftarrow} (x)$ satisfies 
$$
u_{\leftarrow} (x )= \lim_{N\to \infty} e^{-i\kappa N} u_{\leftarrow} (x-Ne_1 )=0.
$$
We can see $ u(x)=0 $ for any $x\in {\bf Z}^2 \setminus \Omega^i $ repeating similar procedures.
\qed

\medskip

The resonances given by Definition \ref{S2_def_resonance} are related to outgoing states as follows.
Let us fix $ \theta \in \mathcal{O}^-_0 $ arbitrarily.
A resonance $e^{-i\kappa} $ with $\kappa\in \mathcal{O}^+ _{\mathrm{Im} \, \theta} $ is characterized by the existence of \textit{outgoing solutions} to the equation $ Uu=e^{-i\kappa} u$, i.e., $u$ is of the form
\begin{gather*}
u_{\leftarrow} (x)= \left\{
\begin{split}
a_{\leftarrow} (x_2) e^{-i\kappa x_1} &, \quad x_1 < -M_0 , \\
0 &, \quad x_1 > M_0 ,
\end{split}
\right. \quad
u_{\rightarrow} (x)= \left\{
\begin{split}
0&, \quad x_1 < -M_0 , \\ 
a_{\rightarrow} (x_2) e^{i\kappa x_1} &, \quad x_1 > M_0 , 
\end{split}
\right.
\end{gather*}
\begin{gather*}
u_{\downarrow} (x)= \left\{
\begin{split}
a_{\downarrow} (x_1) e^{-i\kappa x_2} &, \quad x_2 < -M_0 , \\
0 &, \quad x_2 > M_0 ,
\end{split}
\right. \quad
u_{\uparrow} (x)= \left\{
\begin{split}
0&, \quad x_2 < -M_0 , \\ 
a_{\uparrow} (x_1) e^{i\kappa x_2} &, \quad x_2 > M_0 , 
\end{split}
\right.
\end{gather*}
where $ a_{\leftarrow } (x_2 )$, $a_{\rightarrow} (x_2 )$, $a_{\downarrow} (x_1)$ and $a_{\rightarrow} (x_1 )$ are complex valued sequences.
We show that 
$$ 
\mathrm{supp} \, a_j \subset \{ -M_0 , -M_0 +1 , \ldots ,M_0 -1, M_0 \} , 
$$ 
for every $ j\in \{ \leftarrow , \rightarrow , \downarrow , \uparrow \} $ in the proof of Proposition \ref{S2_lem_outgoing}.
By the definition, outgoing solutions do not belong to $\mathcal{H}$ and they are exponentially growing at infinity if $\kappa\in \mathcal{O}^+ _{\mathrm{Im} \, \theta} \cap \mathcal{O}^-_0$.

\begin{prop}
For $\kappa\in\mathcal{O}^-_0$, $e^{-i\kappa}$ is a resonance of $U$ if and only if there exists a non-trivial outgoing solution $u$ to the equation $ Uu=e^{-i\kappa} u$. 
A solution $u$ is outgoing if and only if $v:=T(\theta)u$ belongs to $\mathcal{H}$ for $\theta\in\mathcal{O}_{\mathrm{Im} \, \kappa}^-$.
\label{S2_lem_outgoing}
\end{prop}

Proof.
If $ e^{-i\kappa} $ is a resonance of $U$, there exists a non-trivial solution $v\in \mathcal{H}$ to the equation $U(\theta)v=e^{-i\kappa} v$ by Definition \ref{S2_def_resonance}.
In view of $U(\theta) = e^{-i\theta} ST(\theta)CT(\theta)^{-1} $, $v$ satisfies $ST(\theta)CT(\theta)^{-1} v= e^{-i(\kappa -\theta)} v$.
For $x_1 < -M_0$, we have 
$$ 
v_{\leftarrow} (x)=e^{-i(\kappa -\theta)} v_{\leftarrow} (x-e_1 ) ,
$$ 
which implies $ v_{\leftarrow} (x) = a_{\leftarrow} (x_2 ) e^{-i (\kappa - \theta )x_1 } $ for a sequence $ a_{\leftarrow} (x_2 ) $. 
For $x_1 > M_0 $, we have 
$$
v_{\leftarrow} ( x+e_1 )= e^{-i(\kappa - \theta )} v_{\leftarrow} (x),
$$
which implies $ v_{\leftarrow} (x)= a'_{\leftarrow} (x_2) e^{-i(\kappa-\theta )x_1} $ for a sequence $ a'_{\leftarrow} (x_2 ) $.
In view of $ v\in \mathcal{H}$, $a' _{\leftarrow} (x_2)$ must vanish for any $x_2$, since $e^{-i (\kappa - \theta) x_1 }$ grows exponentially as $x_1 \to \infty $.
For $ v_j $ with $ j\in \{ \rightarrow , \downarrow , \uparrow \} $, the argument is similar.
We obtain 
\begin{gather}
v_{\leftarrow} (x)= \left\{
\begin{split}
a_{\leftarrow} (x_2 ) e^{-i (\kappa - \theta ) x_1 } &, \quad x_1 < -M_0 , \\
0 &, \quad x_1 > M_0 , 
\end{split}
\right. \label{S2_eq_outgoingleft} \\
v_{\rightarrow} (x)= \left\{
\begin{split}
0 &, \quad x_1 < -M_0 , \\
a_{\rightarrow} (x_2 ) e^{i (\kappa - \theta ) x_1 } &, \quad x_1 >M_0 , \\
\end{split}
\right. \label{S2_eq_outgoingright} 
\end{gather}
\begin{gather}
v_{\downarrow} (x)= \left\{
\begin{split}
a_{\downarrow} (x_1) e^{-i (\kappa - \theta ) x_2 } &, \quad x_2 < -M_0 , \\
0 &, \quad x_2 > M_0 , 
\end{split}
\right. \label{S2_eq_outgoingdown} \\
v_{\uparrow} (x)= \left\{
\begin{split}
0 &, \quad x_2 < -M_0 , \\
a_{\uparrow} (x_1 ) e^{i (\kappa - \theta ) x_2 } &, \quad x_2 >M_0 . \\
\end{split}
\right. \label{S2_eq_outgoingright} 
\end{gather}
For $|x_2| >M_0$, $ v_{\leftarrow} (x+e_1 )= e^{-i (\kappa - \theta )} v_{\leftarrow} (x)$ for all $x_1 \in {\bf Z} $.
Then $a_{\leftarrow} (x_2)=0$ for $|x_2| > M_0 $ due to $v\in \mathcal{H}$.
For $a_{\rightarrow} (x_2)$, $a_{\downarrow} (x_1 )$ and $a_{\uparrow} (x_1)$, the proof is similar.
Now we put $ u=T(\theta)^{-1} v$.
In view of $U(\theta)v=e^{-i\kappa}v$, we have $Uu=e^{-i\kappa} u$.
Due to (\ref{S2_eq_outgoingleft})-(\ref{S2_eq_outgoingright}), $u$ is obviously outgoing.

Let us turn to the proof of the converse.
For an outgoing solution $u$ to $Uu=e^{-i\kappa} u$, we put $v= T(\theta )u$.
Note that we see $\mathrm{supp} \, a_j \subset \{ -M_0 , -M_0 +1 , \ldots ,M_0 -1, M_0 \} $ for every $ j\in \{ \leftarrow , \rightarrow , \downarrow , \uparrow \} $ by the equation $Uu=e^{-i\kappa} u$ and the assumption $\kappa \in \mathcal{O}^+ _{\mathrm{Im} \, \theta} \cap \mathcal{O}^- _0 $.
Due to the definition of $ T(\theta )$, we see that $v\in \mathcal{H}$ solves the equation $ U(\theta)v= e^{-i\kappa} v $.
\qed

\section{Trapping and non-trapping trajectory}
\subsection{Elastic scattering of QWs}
Here we consider a special case of QWs.
Namely, we introduce elastic scattering of QWs and define a trajectory of the quantum walker as a particle.
Let $ \mathcal{P} _{ch} $ be the set of permutations
$$
\sigma = \begin{pmatrix} \leftarrow & \rightarrow & \downarrow & \uparrow \\ \sigma (\leftarrow) & \sigma (\rightarrow) & \sigma (\downarrow) & \sigma (\uparrow) \end{pmatrix}.
$$
We define the coin operator $ C_{el} $ of multiplication by the matrix $ C_{el} (x) \in \mathrm{U} (4) $ of the form
\begin{gather*}
C_{el} (x)= \left\{
\begin{split}
[ e^{i\alpha _j (x)} {\bf e} _{\sigma (x,j)} ]_{j\in \{ \leftarrow , \rightarrow , \downarrow , \uparrow \} } &, \quad x\in \Omega^i , \\
I_4 &, \quad x\in {\bf Z}^2 \setminus \Omega^i ,
\end{split}
\right.
\end{gather*}
where $\alpha_j (x)\in {\bf T}$ and 
$$
\sigma (x):= \begin{pmatrix} \leftarrow & \rightarrow & \downarrow & \uparrow \\ \sigma (x,\leftarrow) & \sigma (x,\rightarrow) & \sigma (x,\downarrow) & \sigma (x,\uparrow) \end{pmatrix} \in \mathcal{P}_{ch} .
$$
The scattering process associated with the QW $ U_{el} =SC_{el} $ is an elastic scattering, i.e., the quantum walker behaves like a classical particle and scatters without loss of its energy at every point $x\in {\bf Z}^2 $.
Indeed, we have 
$$
(C_{el}u)(x)= \sum _{j\in \{ \leftarrow , \rightarrow , \downarrow , \uparrow \} } e^{i\alpha_j (x)} u_j (x) {\bf e} _{\sigma (x,j)}.
$$
It follows that
\begin{gather}
(U_{el} u)_{\leftarrow} (x)= (C_{el} u)_{\leftarrow} (x+e_1) = e^{i\alpha_{\sigma^{-1} (x+e_1 ,\leftarrow)} (x+e_1)} u _{\sigma^{-1} (x+e_1,\leftarrow)} (x+e_1) , \label{S3_eq_elastic1} \\
(U_{el} u)_{\rightarrow} (x)= (C_{el} u)_{\rightarrow} (x-e_1) = e^{i\alpha_{\sigma^{-1} (x-e_1 ,\rightarrow)} (x-e_1)} u _{\sigma^{-1} (x-e_1,\rightarrow)} (x-e_1) , \label{S3_eq_elastic2} \\
(U_{el} u)_{\downarrow} (x)= (C_{el} u)_{\downarrow} (x+e_2) = e^{i\alpha_{\sigma^{-1} (x+e_2 ,\downarrow)} (x+e_2)} u _{\sigma^{-1} (x+e_2 ,\downarrow)} (x+e_2) , \label{S3_eq_elastic3} \\
(U_{el} u)_{\uparrow} (x)= (C_{el} u)_{\uparrow} (x-e_2) = e^{i\alpha_{\sigma^{-1} (x-e_2 ,\uparrow)} (x-e_2)} u _{\sigma^{-1} (x-e_2,\uparrow)} (x-e_2) , \label{S3_eq_elastic4}
\end{gather}
for every $x \in {\bf Z}^2 $, extending $\sigma (x)\in \mathcal{P}_{ch}$ to be the identity for $x\in {\bf Z}^2 \setminus \Omega^i$.

Let us introduce the trajectory of the QW $ U_{el} $.
We take the initial state 
$$
f_{ y,j,\alpha } =\{ \delta_{x,y} e^{i\alpha} {\bf e } _j \} _{x\in {\bf Z}^2} \in \mathcal{H} ,
$$ 
for $y\in {\bf Z}^2 $, $j\in \{ \leftarrow , \rightarrow , \downarrow , \uparrow \} $, and $\alpha \in {\bf T} $.
Here $\delta_{x,y}$ denotes the Kronecker delta.
The dynamics of $ U_{el} $ is an elastic scattering. 
$ U^t_{el} f_{ y,j,\alpha }$ for every $t\in {\bf Z}$ is of the form 
$$
(U^t_{el} f_{ y,j,\alpha } )(x)= \delta_{x,q(t,y,j )} e^{i\alpha (t,y,j)} {\bf e} _{p(t,y,j)} ,
$$
where $\alpha (t,y,j)\in {\bf T} $ and
\begin{gather}
q(t,y,j)= \mathrm{supp} \, (U^t_{el} f_{ y,j,\alpha } ), \quad q(0,y,j)=y, \label{S3_eq_trajectory1} \\
p(t,y,j)\in \{ \leftarrow , \rightarrow , \downarrow , \uparrow \} , \quad p(0,y,j)=j . \label{S3_eq_trajectory2} 
\end{gather}
In view of (\ref{S3_eq_trajectory1})-(\ref{S3_eq_trajectory2}), the mapping $ \Phi ( \cdot ,y,j) : {\bf Z} \to {\bf Z}^2 \times \{ \leftarrow , \rightarrow , \downarrow , \uparrow \} $ defined by 
$$
\Phi (t,y,j)= (q(t,y,j),p(t,y,j)), \quad \Phi (0,y,j)=(y,j),
$$
is one of analogues of classical trajectories.
By the definition, $\Phi (\cdot ,y,j) $ is independent of $\alpha$.
Obviously, the trajectory $\Phi (\cdot,y,j)$ is uniquely determined by the operator $U_{el}$ and the initial value $ (y,j)$.
Thus it follows that each trajectory $ \Phi (\cdot,y,j)$ does not have junctions ($\Phi (\cdot,y,j)$ has neither confluences nor branches).

\begin{definition}
For $ \Phi ( \cdot ,y,j) $, we define the following notions.
\begin{enumerate}
\item We call $ \Phi ( \cdot ,y,j) $ a \textit{bounded trajectory} if there exists a constant $\rho >0$ such that $|q(t,y,j)| \leq \rho $ for any $t\in {\bf Z} $. 

\item We call $ \Phi ( \cdot ,y,j) $ a \textit{closed trajectory} if there exist integers $t_1 , t_2 \in {\bf Z}$ with $ t_1 < t_2 $ such that $ \Phi (t_1 ,y,j) = \Phi (t_2 ,y,j) $.
\end{enumerate}
\label{S3_def_boundedtrajectory}
\end{definition}

For a constant $\rho >0$, the subset $\{ x\in {\bf Z}^2 \ ; \ |x|\leq \rho \} \times \{ \leftarrow , \rightarrow , \downarrow , \uparrow \}  $ is a finite set.
Thus we see the following property.

\begin{lemma}
Let $\Phi (\cdot ,y,j)$ for $(y,j)\in {\bf Z}^2 \times \{ \leftarrow , \rightarrow , \downarrow , \uparrow \} $ be a trajectory.
\begin{enumerate}
\item
The trajectory $\Phi (\cdot ,y,j) $ is closed if and only if $\Phi (\cdot ,y,j) $ is bounded.

\item
The trajectory $\Phi (\cdot ,y,j)$ satisfies $|q(t,y,j)|\to \infty$ as $t\to \infty$ if and only if $|q(t,y,j)|\to \infty$ as $t\to -\infty$.
\end{enumerate}
\label{S3_lem_bounded_closed}
\end{lemma}

Proof.
The statement (1) is trivial in view of the finiteness of the subset $\{ x\in {\bf Z}^2 \ ; \ |x|\leq \rho \} \times \{ \leftarrow , \rightarrow , \downarrow , \uparrow \}  $ for any constant $\rho >0$.
We consider the statement (2).
Suppose that there exists a constant $\rho_+ >0 $ such that $|q(t,y,j)| \leq \rho_+ $ for $t\geq 0$ even though $|q(t,y,j)|\to \infty$ as $t\to -\infty $.
The finiteness of the subset $\{ x\in {\bf Z}^2 \ ; \ |x|\leq \rho_+ \} \times \{ \leftarrow , \rightarrow , \downarrow , \uparrow \}  $ implies that $\Phi (t,y,j)$ for $t\geq 0$ consists of a closed trajectory.
Then we can take nonnegative integers $ t_1 $ and $t_2$ with $t_1 < t_2 $ such that $\Phi (t_1 ,y,j)=\Phi (t_2 ,y,j)$.
Moreover, we can choose $t_1 \geq 0$ as the minimum of integers which satisfy the above situation.
By the assumption, $\Phi (t,y,j)$ for $t< t_1 $ is an unbounded trajectory.
Due to the choice of $t_1$, we can assume that $\Phi (t_1 -1,y,j)\not= \Phi (t_2 -1,y,j)$ even though $\Phi (t_1,y,j)=\Phi (t_2,y,j)$. 
This is a contradiction in view of the definition of $ U_{el} $.
If we replace $t\to \infty$ and $t\to -\infty$, the proof is parallel.
\qed

\medskip

On a closed trajectory $\Phi (\cdot,y,j) $, we can construct an eigenfunction of $U_{el} $ as follows.
The existence of closed trajectories associated with $ U_{el} $ determines the existence of eigenvalues of $ U_{el} $.

\begin{lemma}
There exists a pair $(y,j) \in {\bf Z}^2 \times \{ \leftarrow , \rightarrow , \downarrow , \uparrow \} $ such that $ \Phi ( \cdot , y,j )$ is closed if and only if there exist some eigenvalues of $ U_{el} $. 
If $ \Phi ( \cdot , y,j )$ is closed with its period $N$, then there exist $N$ eigenvalues satisfying \eqref{QC}, where associated eigenfunctions are supported only on $\{\Phi(t,y,j);\,t\ge0\}$. 
\label{S3_lem_evUel}
\end{lemma}

\textit{Remark.} 
Here, we say a sequence $u=\{[u_{\leftarrow} (x),u_{\rightarrow} (x),u_{\downarrow} (x),u_{\uparrow} (x)]^{\mathsf{T} } \}_{x\in{\bf Z}^2}$ is supported on a subset $A\subset{\bf Z}^2\times \{ \leftarrow , \rightarrow , \downarrow , \uparrow \} $ if $u_j(x)=0$ for each $(x,j)\notin A$. 
Note that the condition \eqref{QC} associated with the periodic trajectory is similar to the well-known Bohr-Sommerfeld quantization condition in the quantum mechanics.

\medskip

Proof.
Suppose that there is an initial value $(y,j)$ such that $\Phi (\cdot , y,j)$ is closed.
Without loss of generality, we can assume that $ \Phi (0,y,j)=\Phi (N,y,j)$ for a positive integer $ N$.
In fact, $N$ must be even.
In the following, we choose the smallest $N$ such that the above situation holds. 
Let $ q(t) = q(t,y,j) $ and $p(t) = p(t,y,j) $ for $t=0,1,2,\ldots$.
Thus the trajectory $ \{ (q(t),p(t)) \} _{t\geq 0}$ is a closed with $(q(0),p(0))=(q(N),p(N))=(y,j)$.
By the definition of $\Phi (\cdot ,y,j) $, we note 
$$
(U_{el} u)_{p(t) } (q(t) )= e^{i\alpha _{p(t-1)} (q(t-1))} u_{p(t-1)} (q(t-1) ), 
$$
for every $ t=0,1,\ldots,N-1$.
Letting $ f_t = u_{p(t)} (q(t) ) $ and $\beta_t = \alpha _{p(t)} (q(t) )$, the equation $U_{el} u= e^{-i\lambda} u$ on $\{ q(t) \} _{t=0}^N $ can be rewritten as 
\begin{equation}
e^{i\beta _{t-1}} f_{t-1} = e^{-i\lambda} f_t ,
\label{S3_eq_evUel01}
\end{equation}
in view of (\ref{S3_eq_elastic1})-(\ref{S3_eq_elastic4}).
As a consequence, we have the equality
$$
f_0 = e^{i(\lambda N+\beta_0 + \cdots + \beta_{N-1})} f_0 .
$$
Suppose that $ f_0 \not= 0$ is given and the equality 
\begin{equation}
\lambda \equiv - \frac{1}{N} \sum _{t=0}^{N-1} \beta_t \quad \text{modulo} \quad 2\pi ,
\label{QC}
\end{equation}
holds true.
Then the sequence $\{ f_t \} $ is uniquely determined from $f_0$ by the equality (\ref{S3_eq_evUel01}).
Moreover, the function $u\in \mathcal{H}$ with 
\begin{gather}
u_j (x)= \left\{
\begin{split}
f_t &, \quad x=q(t) , \quad j=p(t) , \quad t=0,1,\ldots,N-1, \\
0 &, \quad \text{otherwise} , 
\end{split}
\right.
\label{EFalongCT}
\end{gather}
is an eigenfunction of $ U_{el} $ associated with the eigenvalue $e ^{-i\lambda} $.

Conversely, suppose that for any $(y,j)$, the trajectory $\Phi (\cdot ,y,j) $ is not closed.
Let $u\in \mathcal{H}$ be a solution to the equation $U_{el} u= e^{-i\lambda} u $ for a constant $\lambda \in {\bf T} $. We prove the argument by showing that $u$ is trivial, that is, $u$ vanishes identically.
Fix an arbitrary point $x\in {\bf Z}^2$ and an arbitrary chirality $j\in \{ \leftarrow , \rightarrow , \downarrow , \uparrow \} $. 
We consider the trajectory $\Phi (\cdot ,x,j) $.
Since there is no bounded trajectory in view of Lemma \ref{S3_lem_bounded_closed}, we note $ |q(t,x,j)| \to \infty$ as $|t|\to \infty $.
We introduce the notation $\Phi (t,x,j) =  (q(t),p(t)) $ by $q(0) =x$, $q(\pm 1) =q ( \pm 1,x,j)$, $q(\pm 2) = q(\pm 2,x,j)$, $\ldots $ and $p(0)=j$, $p(\pm 1) =p(\pm 1,x,j)$, $p(\pm 2) = p(\pm 2 ,x,j)$, $\ldots$.
Letting $ f_t = u_{p(t)} (q(t) ) $ and $\beta_t = \alpha _{p(t)} (q(t) )$, we obtain 
\begin{equation}
e^{ i\beta_{t-1} } f_{t-1} = e^{-i\lambda} f_t , \quad t\in {\bf Z},
\label{S3_eq_evUel03}
\end{equation}
and thus 
\begin{equation}
f_0 = e^{-i (\lambda N + \beta_0 + \beta_1 +\cdots + \beta_{N-1} )} f_N,\quad
\left|f_N\right|=\left|f_0\right|
\label{S3_eq_evUel02}
\end{equation}
for any positive integer $N$. 
On the other hand, $u\in\mathcal{H}$ implies 
\begin{equation*}
\sum_{t\in{\bf Z}} \left| f_0 \right|^2  = \sum_{t\in{\bf Z}}\left|f_t\right|^2\leq \| u \|^2 _{\mathcal{H}} < + \infty .
\end{equation*}
As a consequence, we obtain $f_0 =u_j (x)=0$. We conclude that $u=0$ since the choice of $(x,j)$ is arbitrary.
\qed

\medskip

In view of Lemma \ref{S3_lem_evUel}, we call the QW $U_{el} $ \textit{non-trapping} if there is no closed trajectory associated with $U_{el} $.
We saw that a non-trapping elastic QW has no eigenvalue. The following lemma shows that it has no resonance. Moreover, the absence of resonances (other than eigenvalues) of elastic QWs is always true.

\begin{lemma}
Elastic QWs have no resonance other than eigenvalues. Moreover, a non-trapping QW does not have eigenvalues either.
\label{S3_lem_NTresonance}
\end{lemma}

Proof. 
Let an outgoing sequence $u$ satisfy $e^{-i\kappa}u=U_{el}u$ for some $\kappa\in\mathcal{O}_0^-$. For each $(x,j)\in{\bf Z}^2\times\{\leftarrow,\rightarrow,\downarrow,\uparrow\}$ such that $\Phi(\cdot,x,j)$ is unbounded, $u_j(x)=0$ follows.  
In particular, the absence of resonances of non-trapping elastic QWs is already obtained. Here, we use the facts that $u$ is outgoing and that $\Phi(\cdot,x,j)$ is an incoming straight line for $t\in(-\infty,t_0]$ for some $t_0\in{\bf Z}$. 

If $u_j(x)\neq0$ holds for some $(x,j)$ such that $\Phi(\cdot,x,j)$ is bounded, the same argument as the proof of Lemma \ref{S3_lem_evUel} shows that $\kappa$ has to satisfy the quantization condition \eqref{QC}. This implies that $\kappa$ is real, and contradicts with $\kappa\in\mathcal{O}_0^-$.
\qed

\subsection{Eigenvalue of QWs via non-penetrable barriers}
Let us consider another setting of QWs which have no resonances other than some eigenvalues.
In our model, we can know the sum of geometric multiplicities of eigenvalues.
Namely, we introduce a ``non-penetrable barrier" on the boundary of $ \Omega^i $.
Let $U_{np}$ denote a QW satisfying
\begin{gather}
\begin{split}
&(U u)_{\rightarrow } (x+e_1 )= u_{\leftarrow } (x) , \quad x\in K_1^+ ,  \\ 
&(U u)_{\leftarrow } (x-e_1 )= u_{\leftarrow } (x) , \quad x\in K_1^- ,  \\
&(U u)_{\uparrow } (x+e_2 )= u_{\downarrow } (x) , \quad x\in K_2^+ ,  \\ 
&(U u)_{\downarrow } (x-e_2 )= u_{\uparrow } (x) , \quad x\in K_2^- , 
\end{split}
\label{S3_eq_np3}
\end{gather}
in addition to the assumption (A-1), where we put
\begin{gather*}
K_1^{\pm} = \{ x\in \Omega^i \ ; \ x_1 = \pm M_0 , \ -M_0 \leq x_2 \leq M_0  \}  , \\
K_2^{\pm} =  \{ x\in \Omega^i \ ; \ x_2 = \pm M_0 , \ -M_0 \leq x_1 \leq M_0  \}  .
\end{gather*}
Under above boundary condition on $K= K^+_1 \cup K^- _1 \cup K^+ _2 \cup K_2^- ,$  $U_{np}$  is completely split into two QWs $U_i$ and $U_e$ independent of each other in the following manner. Note that $ (\pm M_0 , \pm M_0 ) \in K_1^{\pm} \cap K_2 ^{\pm} $ and $(\pm M_0 , \mp M_0 )\in K_1^{\pm} \cap K_2^{\mp} $, respectively.

\begin{figure}[b]
\centering
\includegraphics[bb=0 0 659 655, width=6cm]{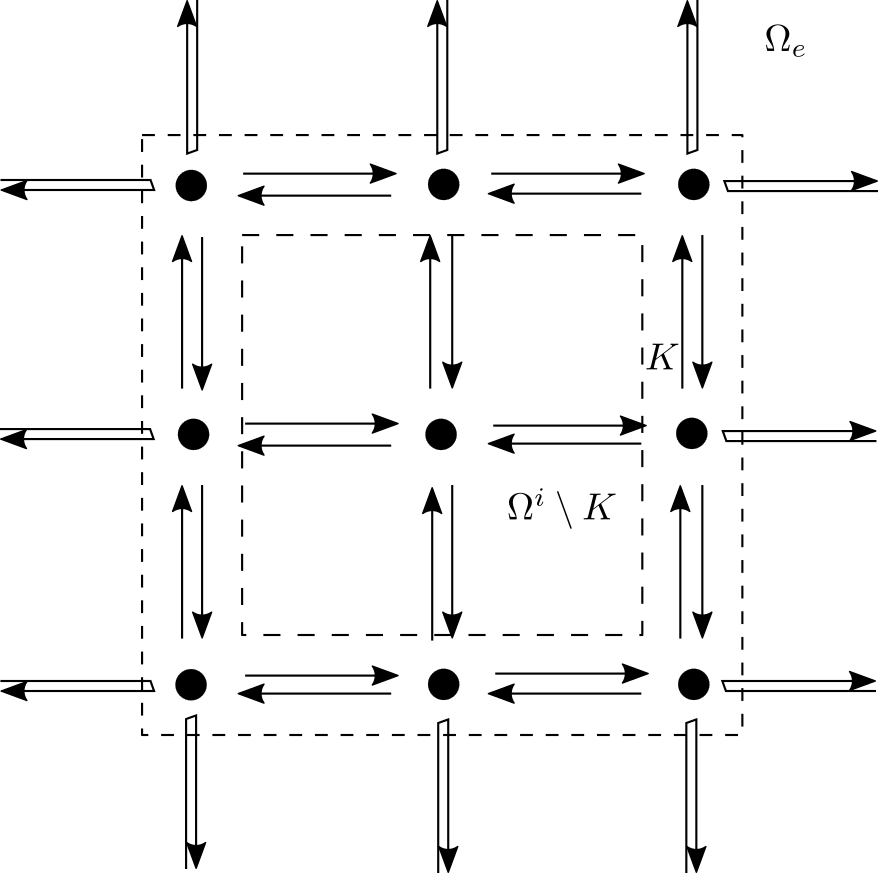}
\caption{The situation of $\Omega^i$, $\Omega^e$, and $K$ for $M_0 = 1$.}
\label{fig_QW_Unp}
\end{figure}

We define the exterior domain $ \Omega^e $ by
$$
\Omega^e = ( {\bf Z}^2 \setminus  \Omega^i  ) \cup K.
$$
Now we decompose $\mathcal{H}$ into the direct sum of Hilbert spaces $\mathcal{H}_i \oplus \mathcal{H}_e $ where
\begin{gather*}
\mathcal{H}_i = \{ u\in \ell^2 ( {\bf Z}^2 ; {\bf C}^4 ) \ ; \ u\text{ satisfies the condition } (\ref{S3_eq_bcint}) \} , \\
\mathcal{H}_e = \{ u\in \ell^2 ( {\bf Z}^2 ; {\bf C}^4 ) \ ; \ u\text{ satisfies the condition } (\ref{S3_eq_bcext}) \}  ,
\end{gather*} 
with
\begin{equation}
\mathrm{supp} \, u \subset \Omega^i, \quad u_{\leftarrow} | _{K_1^+} = u_{\rightarrow} | _{K_1^-} = u_{\downarrow} | _{K_2^+} = u _{\uparrow} | _{K_2^-} =0, 
\label{S3_eq_bcint}
\end{equation}
\begin{gather}
\mathrm{supp} \, u \subset \Omega^e , \quad
\left\{
\begin{split}
&u_{\leftarrow} | _{K_1^- \cup K_2^+ \cup K_2^- } = u_{\rightarrow} | _{K_1^+ \cup K_2^+ \cup K_2^- } =0 ,  \\
&u_{\downarrow} | _{K_2^- \cup K_1^+ \cup K_1^- } = u_{\uparrow} | _{K_2^+ \cup K_1^+ \cup K_1^- } =0.
\end{split}
\right. 
\label{S3_eq_bcext}
\end{gather}
Then we have $ U _{np} = U_i \oplus U_e $ on $\mathcal{H}_i \oplus \mathcal{H}_e $ with $U_i=U|_{\mathcal{H}_i}$ and $U_e=U|_{\mathcal{H}_e}$. 
Roughly speaking, the condition (\ref{S3_eq_np3}) means that the boundary $K$ acts the reflector from $\Omega^e$ to $\Omega^e$ as well as $\Omega^i$ does not leak quantum walkers from $\Omega^i $.
See Figure \ref{fig_QW_Unp}.

\medskip

\textit{Remark.}
The boundary condition defined by (\ref{S3_eq_bcint}) and (\ref{S3_eq_bcext}) is slightly complicated.
In order to avoid the complexity, it is convenient to identify $\Omega^i$ and $\Omega^e $ with graphs $\Gamma^i $ and $\Gamma^e $ and consider the corresponding QWs on $\Gamma^i$ and $\Gamma^e $ as follows.
Let $\Gamma^i = (\mathcal{V}^i , \mathcal{A}^i ) $ be a finite graph where $\mathcal{V}^i = \Omega^i$ is the set of vertices and $\mathcal{A}^i$ is the set of oriented edges $\omega = ( o(\omega),t(\omega)) $ with $o(\omega),t(\omega) \in \Omega^i $ and $|o(\omega)-t(\omega)|=1$.
Here $o(\omega)$ and $t(\omega)$ denote the origin and the terminus of the edge $\omega$, respectively.
The infinite graph $\Gamma^e = (\mathcal{V}^e , \mathcal{A}^e )$ corresponds to the exterior domain $\Omega^e $.
Here $\mathcal{V}^e = \Omega^e$ and $\mathcal{A}^e $ is the set of oriented edges $\omega = (o(\omega),t(\omega))$ such that both of neighboring vertices $o(\omega)$ and $t(\omega)$ belong to $\Omega^e \setminus K$ or one of neighboring vertices $o(\omega)$ and $t(\omega)$ lies in $K$ and the other lies in $\Omega^e \setminus K$.

Let $\ell^2 (\mathcal{A}^i )$ be the Hilbert space equipped with the inner product
$$
(\widetilde{u},\widetilde{v} )_{\ell^2 (\mathcal{A}^i ) }  = \sum _{\omega \in \mathcal{A}^i} \widetilde{u} (\omega) \overline{\widetilde{v} (\omega)} , \quad \widetilde{u},\widetilde{v} \in \ell^2 (\mathcal{A}^i ) .
$$
Now we consider the identification by the unitary transform $\mathcal{I}_i :\mathcal{H}_i \to \ell^2 (\mathcal{A}^i )$ as 
\begin{gather*}
(\mathcal{I}_i u)(x+e_1,x)=u_{\leftarrow} (x), \quad (\mathcal{I}_i u)(x-e_1,x)=u_{\rightarrow} (x), \\
(\mathcal{I}_i u)(x+e_2,x)=u_{\downarrow} (x), \quad (\mathcal{I}_i u)(x-e_2,x)=u_{\uparrow} (x),
\end{gather*}
if $ (x\pm e_1,x), (x\pm e_2,x)\in \mathcal{A}^i$.
Then $\widetilde{U}_i := \mathcal{I}_i U_i \mathcal{I}_i^{-1} $ is a QW on the finite graph $\Gamma^i $ without any boundary condition.
We also see that $\widetilde{U}_i $ is unitary on $\ell^2 (\mathcal{A}^i )$.
For $U_e$ and $\Gamma^e $, the argument is similar.

\medskip

In view of the above remark, we see the following result on the spectrum of $U_i$.

\begin{lemma}
We have $\sigma (U_i) = \sigma_p (U_i) \subset S^1 $.
Each eigenvalue has a finite multiplicity.
The sum of geometric multiplicities of eigenvalues coincides with $ N= \# \mathcal{A}^i $.
Namely, there exist orthonormal eigenfunctions $u^{(1)} , \ldots , u^{(N)} $ in $\mathcal{H}_i$.
\label{S3_lem_evUnp}
\end{lemma}

Proof.
Noting that the graph $\Gamma^i$ is finite and $ \widetilde{U}_i $ is unitary on $\ell^2 (\mathcal{A}^i )$, this lemma is a direct consequence of the above remark.
\qed

\medskip

Let us turn to the representation of the Green operator
$$
R_i (\kappa )= (U_i -e^{-i\kappa} )^{-1} ,
$$
by the eigenvalues and the associated eigenfunctions.
Suppose that $ e^{-i\kappa_1} , \ldots, e^{-i\kappa_N} \in \sigma_p (U_i) $ and take orthonormal eigenfunctions $u^{(1)} , \ldots ,u^{(N)} \in \mathcal{H}_i $.
Letting $u=R_i (\kappa )f$ for $f\in \mathcal{H}_i$ and $e^{-i\kappa} \not\in \sigma_p (U_i)$, we obtain 
\begin{equation}
u= \sum _{j=1}^N \frac{(f,u^{(j)} ) _{\mathcal{H}_i}}{e^{-i\kappa_j} - e^{-i\kappa}} u^{(j)} ,
\label{S3_eq_Greenf}
\end{equation}
by a direct calculation.
For an eigenvalue $e^{ -i\mu }\in \sigma_p (U_i)$, let us show an estimate of $R_i (\kappa )$ near $\mu \in {\bf T}$.
We take a counterclockwise loop $\mathcal{L}_{\epsilon,s} (\mu )=\sum _{j=1}^4 \mathcal{L} _{\epsilon ,s,j} (\mu )$ where 
\begin{gather}
\begin{split}
&\mathcal{L}_{\epsilon,s,1} (\mu )= \{ (1-\tau )(\mu +a\epsilon^s -ib \epsilon^s ) + \tau (\mu +a\epsilon^s +ib \epsilon^s ) \ ; \ \tau \in [0,1] \} , \\
&\mathcal{L}_{\epsilon,s,2} (\mu )= \{ (1-\tau )(\mu +a\epsilon^s +ib \epsilon^s ) + \tau (\mu -a\epsilon^s +ib \epsilon^s ) \ ; \ \tau \in [0,1] \} ,  \\
&\mathcal{L}_{\epsilon,s,3} (\mu )= \{ (1-\tau )(\mu -a\epsilon^s +ib \epsilon^s ) + \tau (\mu -a\epsilon^s -ib \epsilon^s ) \ ; \ \tau \in [0,1] \} ,  \\
&\mathcal{L}_{\epsilon,s,4} (\mu )= \{ (1-\tau )(\mu -a\epsilon^s -ib \epsilon^s ) + \tau (\mu +a\epsilon^s -ib \epsilon^s ) \ ; \ \tau \in [0,1] \} , 
\end{split}
\label{S2_eq_loop01}
\end{gather}
for some constants $a,b,s,\epsilon >0$.
Taking sufficiently small $\epsilon >0 $, we can assume that there is no other $\kappa$ inside $\mathcal{L}_{\epsilon,s} (\mu )$ such that $e^{-i\kappa} \in \sigma_p (U_i)$.

\begin{lemma}
Let $e^{-i\mu} \in \sigma_p (U_i)$.
If $\kappa \in {\bf T} _{{\bf C}} $ varies on the loop $\mathcal{L} _{\epsilon,s} (\mu )$, we have $\| R_i (\kappa) \| _{{\bf B} (\mathcal{H}_i )} = O(\epsilon^{-s} )$ as $\epsilon \downarrow 0$.
\label{S3_lem_estimateRi}
\end{lemma}

Proof.
We have 
$$
| e^{-i\mu} - e^{-i\kappa} |^2 = 4e^{\mathrm{Im} \, \kappa} |\sin ((\kappa- \mu)/2)|^2 .
$$
Thus, for sufficiently small $|\kappa-\mu|$, there exists a constant $\delta >0$ such that 
\begin{equation}
| e^{-i\mu} - e^{-i\kappa} | \geq \delta |\kappa -\mu | .
\label{S3_eq_lbound_e}
\end{equation}
Suppose that $\kappa \in \mathcal{L} _{\epsilon,s} (\mu )$.
For the case $ \kappa \in \mathcal{L} _{\epsilon,s,1} (\mu ) \cup \mathcal{L} _{\epsilon,s,3} (\mu )$, we have $ | \mathrm{Re} \, \kappa - \mu | =a\epsilon^s $ and $|\mathrm{Im} \, \kappa | \leq b \epsilon^s $.
If $ \kappa \in \mathcal{L} _{\epsilon,s,2} (\mu ) \cup \mathcal{L} _{\epsilon,s,4} (\mu )$, we have $| \mathrm{Im} \, \kappa | = b\epsilon^s  $ and $ |\mathrm{Re} \, \kappa |\leq a \epsilon^s $.
Then, by the formulas (\ref{S3_eq_Greenf}) and (\ref{S3_eq_lbound_e}), we can take a constant $\delta >0$ such that $ \| R_i (\kappa ) f \| _{\mathcal{H}_i} \leq \delta \epsilon^{-s} \| f\| _{\mathcal{H}_i} $ for any $f\in \mathcal{H}_i$.
\qed

\medskip

For $f\in \mathcal{H}_i$, the complex translation $T(\theta)f \in \mathcal{H}_i$ can be defined naturally in view of (\ref{S2_eq_complextranslation_lattice}).
Then we consider the operator $ U_i (\theta)=T(\theta)U_i T(\theta)^{-1} $ and $R_i (\kappa ,\theta)= (U_i (\theta) -e^{-i\kappa} )^{-1} $ for some $\kappa \in {\bf T} _{{\bf C}} $.

\begin{lemma}
Fix $ \theta \in \mathcal{O}^-_0$.
\begin{enumerate}
\item For each $\kappa\in{\bf T}_{\bf C}$ with $\mathrm{Im}\,\kappa\ge0$, $e^{-i\kappa} \in \sigma_p (U_i) $ if and only if $e^{-i\kappa} \in \sigma_p (U_i (\theta)) $. In particular, $U_i$ has no resonances other than eigenvalues.

\item For $f\in \mathcal{H}_i$ and $\kappa\in\mathcal{O}^+_{\mathrm{Im}\,\theta}$ with $e^{-i\kappa} \not\in \sigma_p ( U_i (\theta )) $, we have 
$$
R_i (\kappa , \theta)f = \sum _{j=1}^N \frac{(T(\theta)^{-1} f,u^{(j)} ) _{\mathcal{H}_i}}{e^{-i\kappa_j} - e^{-i\kappa}} T(\theta)u^{(j)}.
$$
Suppose that $\epsilon >0$ is sufficiently small.
If $\kappa$ varies on $\mathcal{L} _{\epsilon,s} (\mu )$ for an eigenvalue $e^{-i\mu} \in \sigma_p (U_i(\theta))$, we have $ \| R_i (\kappa , \theta) \| _{{\bf B} (\mathcal{H}_i)} =O(\epsilon^{-s})$.

\end{enumerate}
\label{S3_lem_complexUi}
\end{lemma}

Proof.
The assertion (1) follows from
$$
U_i u=e^{-i\kappa} u \iff U_i (\theta)v=e^{-i\kappa}v,
$$
for $v= T(\theta)u$. According to Lemma \ref{S3_lem_evUnp}, the eigenfunctions of $U_i$ spans $\mathcal{H}_i$. Thus, there is no resonance of $U_i$ other than eigenvalues.

Suppose that $u,f\in \mathcal{H}_i$ satisfy $(U_i (\theta) -e^{-i\kappa} )u=f$.
This equation is equivalent to $(U_i -e^{-i\kappa} ) T(\theta)^{-1} u=T(\theta )^{-1} f$.
We apply the formula (\ref{S3_eq_Greenf}) to this equation.
The representation in the assertion (2) follows.
The remaining part is parallel to Lemma \ref{S3_lem_estimateRi}.
\qed

\medskip

Let us turn to the exterior QW $U_e$ and its resolvent operator
$$
R_e (\kappa) = (U_e -e^{-i\kappa} )^{-1} .
$$
We define $ U_e (\theta)= T(\theta)U_e T(\theta)^{-1} $ and $ R_e (\kappa , \theta)= (U_e ( \theta)-e^{-i\kappa} )^{-1} $ for some $ \kappa \in {\bf T} _{{\bf C}} $.
The following estimate is used in Section 5.

\begin{lemma}
Fix $\theta \in \mathcal{O}^-_0$ and a compact subset $Z \subset \mathcal{O}^+ _{\mathrm{Im} \, \theta}$.
Then there exists a constant $\gamma >0$ which depends only on $Z$ such that
$$
\| R_e (\kappa , \theta) \| _{{\bf B} (\mathcal{H}_e)} \leq \gamma ,
$$
for $\kappa \in Z$.
As a consequence, there is no resonance of $U_e$.
\label{S3_lem_estimateRe}
\end{lemma}

Proof.
Let us introduce the operator of restriction $ J_e : \mathcal{H} \to \mathcal{H}_e $.
Thus $J_e$ is the orthogonal projection onto $\mathcal{H}_e$.
Its adjoint $J_e^* $ satisfies $ J^*_e : \mathcal{H}_e \ni g \mapsto g\in \mathcal{H} $.
Letting $ u= R_e ( \kappa , \theta )f$ for $f\in \mathcal{H}_e$, we have 
$$
(U_0 (\theta) -e^{-i\kappa} ) J_e^* u =f+Q(\theta)u,
$$
where $Q(\theta)= U_0 (\theta) J_e^* - J_e^* U_e (\theta )$.
It follows from this equality that 
\begin{equation}
J_e^* R_e (\kappa, \theta )= R_0 (\kappa , \theta ) + R_0 (\kappa , \theta) Q(\theta) R_e (\kappa , \theta) .
\label{S3_eq_resolventeqext}
\end{equation}
If the estimate of the lemma fails, there exist sequences $\{ f^{(j)} \} _{j\geq 1} \subset \mathcal{H}_e $, $\{ \kappa_j \} _{j\geq 1} \subset Z$, and $\kappa \in Z$ such that $\| f^{(j)} \| _{\mathcal{H}_e} \to 0$, $\| R_e (\kappa _j , \theta )f^{(j)} \| _{\mathcal{H}_e} =1$, and $\kappa_j \to \kappa $ as $j\to \infty$.
We put $u^{(j)} = R_e ( \kappa_j , \theta )f^{(j)} $.
Since the operator $Q(\theta)$ is of finite rank, there exists a subsequence $ \{ u^{(j_k)} \} _{k\geq 1} $ such that $Q(\theta) u^{(j_k)} \to g$ as $k\to \infty$ for a function $g\in \mathcal{H}$.
Obviously, $U_0$ is non-trapping.
Then Lemma \ref{S3_lem_NTresonance} implies $ R_0 (\kappa _{j_k} ,\theta ) \to R_0 (\kappa , \theta )$ in ${\bf B} (\mathcal{H})$ as $k\to \infty $.
Then it follows from (\ref{S3_eq_resolventeqext}) that
$$
v:= \lim _{k\to \infty} J_e^* u^{(j_k)} = R_0 (\kappa , \theta) g \quad \text{in} \quad \mathcal{H} .
$$
We put $u_e =J_e v = \lim_{k\to \infty} u^{(j_k)} \in \mathcal{H}_e$.
By the definition of $u^{(j_k)}$, we have $ (U_e (\theta) -e^{-i\kappa} ) u_e = 0 $.
However, $U_e$ is also non-trapping since $U_e$ does not have bounded trajectories in $\Omega^e$.
Then we see $u_e =0$ by the same way of the proof of Lemma \ref{S3_lem_NTresonance}.
This fact contradicts with $\| u^{(j_k)} \| _{\mathcal{H}_e} =1$ for any $k$.
\qed

\section{Resonance associated with perturbed closed trajectories}
\subsection{Perturbation of closed trajectories}

\begin{figure}[b]
\centering
\includegraphics[bb=0 0 659 653, width=6cm]{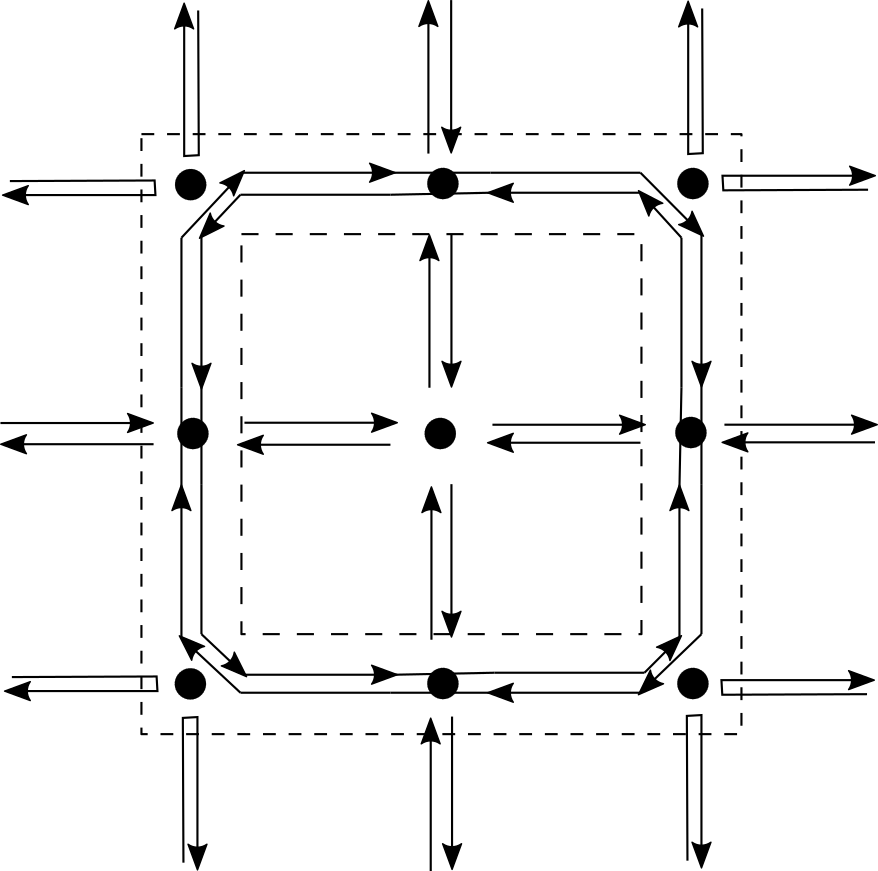}
\caption{The situation of $U_{el}$ for the case $m_0 =n_0 =2 $.}
\label{fig_Uelexample}
\end{figure}

In this section, an example of resonances is given in a constructive approach for the case where a QW $U$ is a small perturbation of $U_{el}$ with only a few closed trajectories.
For the sake of simplicity, we consider a simple model as follows.
Fix four points $(0,0)$, $(m_0, 0)$, $(m_0 ,n_0 )$, $(0,n_0 ) $ for two positive integers $m_0$ and $n_0 $.
We introduce $C _{el} (x)$ for $x\in {\bf Z}^2$ by
\begin{gather}
C_{el} (x)= \left\{
\begin{split}
I_4 &, \quad x\in {\bf Z}^2 \setminus \{ (0,0), (m_0, 0), (m_0 ,n_0 ), (0,n_0 )  \} , \\
[ {\bf e} _{\uparrow} , {\bf e} _{\leftarrow} , {\bf e} _{\rightarrow} , {\bf e} _{\downarrow} ] &, \quad x= (0,0), \\
[ {\bf e} _{\rightarrow} , {\bf e} _{\uparrow} , {\bf e} _{\leftarrow} , {\bf e} _{\downarrow} ] &, \quad x= (m_0 ,0), \\
[ {\bf e} _{\rightarrow} , {\bf e} _{\downarrow} , {\bf e} _{\uparrow} , {\bf e} _{\leftarrow} ] &, \quad x= (m_0 ,n_ 0), \\
[ {\bf e} _{\downarrow} , {\bf e} _{\leftarrow} , {\bf e} _{\uparrow} , {\bf e} _{\rightarrow} ] &, \quad x= (0,n_0 ).
\end{split}
\right.
\label{S4_eq_Celexample}
\end{gather}
See Figure \ref{fig_Uelexample}.
Then $U_{el}$ satisfies (A-1) and has only two closed trajectories $\Phi_+ (\cdot)=(q_+(\cdot),p_+(\cdot)):=\Phi(\cdot ,0,\leftarrow )$ and $\Phi_-(\cdot)=(q_- (\cdot ),p_- (\cdot)):=\Phi (\cdot,0,\downarrow)$ where
\begin{gather*}
q_+ ( t)= \left\{
\begin{aligned}
&(0,t), && \quad 0\leq t\leq n_0 , \\ 
&(t-n_0,n_0 ), &&  \quad n_0 <t\leq n_0 +m_0 , \\ 
&(m_0 , 2n_0+m_0 -t), && \quad n_0 +m_0 < t\leq 2n_0 + m_0 , \\
&(2(n_0+m_0) -t,0), && \quad 2n_0 +m_0<t\leq 2(n_0 +m_0 ) ,
\end{aligned}
\right.
\end{gather*}
\begin{gather*}
p_+ ( t)= \left\{
\begin{split}
\uparrow &, \quad 0<t\leq n_0 , \\ 
\rightarrow & , \quad n_0 <t\leq n_0 +m_0 , \\ 
\downarrow &, \quad n_0 +m_0 < t\leq 2n_0 + m_0 , \\
\leftarrow &, \quad 2n_0 +m_0<t\leq 2(n_0 +m_0 ) ,
\end{split}
\right.
\end{gather*}
and $\Phi_- (\cdot)=\Phi(\cdot ,0,\downarrow )$ is the inverse path of $\Phi _+ (\cdot)=\Phi(\cdot ,0,\leftarrow ) $.

As we have seen in Lemma \ref{S3_lem_evUel}, the set of eigenvalues is given by
\begin{equation}
\sigma_p(U_{el})= \{ e^{-i\pi k/ (m_0 +n_0 )} \ ; \ k= 0,1,\ldots, 2(m_0 +n_0 )-1 \} .
\label{EV-Uel}
\end{equation}
Each eigenvalue has the geometric multiplicity $2$ corresponding to the number of the periodic paths which have the same period as each other. 

\medskip

Now we take a family of QWs $\{U_{\epsilon};\epsilon\in(0,1]\}$ with $U_{\epsilon}=SC_{\epsilon} $ such that each matrix $ C_{\epsilon} (x) \in \mathrm{U} (4)$ satisfies 
\begin{gather}
C_{\epsilon} (x)=C_{el}(x)=I_4 , \quad x\in {\bf Z}^2 \setminus \{ (0,0), (m_0, 0), (m_0 ,n_0 ), (0,n_0 )  \} , \label{S4_eq_Ueassumption1} \\
|C_{\epsilon} (x)- C_{el} (x)| _{\infty} < \epsilon , \quad x\in \{ (0,0), (m_0, 0), (m_0 ,n_0 ), (0,n_0 )  \},\label{S4_eq_Ueassumption2}\\
\begin{aligned}
&c_{\epsilon,\rightarrow,\leftarrow}(0,0)=c_{\epsilon,\uparrow,\downarrow}(0,0)
=c_{\epsilon,\uparrow,\downarrow}(m_0,0)=c_{\epsilon,\leftarrow,\rightarrow}(m_0,0)\\
&=c_{\epsilon,\downarrow,\uparrow}(m_0,n_0)=c_{\epsilon,\leftarrow,\rightarrow}(m_0,n_0)
=c_{\epsilon,\rightarrow,\leftarrow}(0,n_0)=c_{\epsilon,\downarrow,\uparrow}(0,n_0)=0,
\end{aligned}\label{S4_eq_Ueassumption3}
\end{gather}
for each $\epsilon\in(0,1]$, where $c_{\epsilon,j,k}(x)$ $(j,k \in \{ \leftarrow , \rightarrow , \downarrow , \uparrow \})$ denotes the $(j,k)$-entry of $C_{\epsilon} (x)$.
Here we used the norm $|A|_{\infty} = \max _{j,k} |a_{j,k} |$ for a matrix $A=[a_{j,k} ]_{j,k \in \{ \leftarrow , \rightarrow , \downarrow , \uparrow \} } $.

\subsection{Construction of outgoing solutions}
We here discuss the eigenvalues and the resonances of the QW $U_\epsilon$ for each fixed $\epsilon$. The following proposition gives the quantization conditions along $\Phi_+$ of eigenvalues and resonances, which is an analogue of \eqref{QC}.  

\begin{prop}
For each $\epsilon$, there exist eigenvalues of $U_\epsilon$ whose associated eigenfunctions are supported only on the image of $\Phi_+$ if and only if 
\begin{equation}
\left|c_\epsilon^+\right|=1,\quad
c_\epsilon^+:=
c_{\epsilon,\uparrow,\leftarrow} (0,0) c_{\epsilon,\leftarrow,\downarrow} (m_0 ,0)  c_{\epsilon,\downarrow,\rightarrow} (m_0 ,n_0 )c_{\epsilon,\rightarrow,\uparrow} (0,n_0 ).
\label{EVorRes}
\end{equation}
Under this condition, such eigenvalues are  $2(m_0+n_0)$-roots belonging to ${\bf T}$ of
\begin{equation}
e^{-2i(m_0+n_0)\kappa}=c_\epsilon^+.
\label{QC2}
\end{equation}
Otherwise, $2(m_0+n_0)$-roots belonging to $\mathcal{O}_0^-$ of \eqref{QC2} 
are resonances. Associated resonant states are supported on the union of the image of $\Phi_+$ and the eight outgoing tails starting from the four corners $(0,0)$, $(m_0,0)$, $(m_0,n_0)$, $(0,n_0)$, i.e.,
\begin{align*}
&\{((-N,y_2),\leftarrow)\}_{N\ge1},\ 
&&\{((y_1,n_0+N),\uparrow)\}_{N\ge1},\\
&\{((y_1,-N),\downarrow)\}_{N\ge1},\ 
&&\{((m_0+N,y_2),\rightarrow)\}_{N\ge1},
\end{align*}
for $y_1\in\{0,m_0\}$, $y_2\in\{0,n_0\}$.
\label{EVorRes-prop}
\end{prop}

Proof.
The eigenfunctions are constructed in the same way as Lemma~\ref{S3_lem_evUel}. Let us suppose $U_{el}u=e^{-i\kappa}u$ for a $\kappa\in{\bf T}_{\bf C}$. We put $f_t=u_{p_+(t)}(q_+(t))$ for $(q_+(t),p_+(t))=\Phi_+(t)$ defined by \eqref{S4_eq_Celexample}. Then we have $f_t=e^{-i\kappa}f_{t+1}$ for $t\notin\{0,n_0,m_0+n_0,m_0+2n_0\}$ and 
\begin{gather*}
\begin{aligned}
&c_{\epsilon,\uparrow,\leftarrow} (0,0) f_0 = e^{-i\kappa} f_1 , 
&&c_{\epsilon,\rightarrow,\uparrow} (0,n_0 )f_{n_0}= e^{-i\kappa} f_{n_0+1}, \\
&c_{\epsilon,\downarrow,\rightarrow} (m_0 ,n_0 ) f_{m_0 +n_0} = e^{-i\kappa} f_{m_0 +n_0 +1} , &&c_{\epsilon,\leftarrow,\downarrow} (m_0 ,0 ) f_{m_0 +2n_0} = e^{-i\kappa} f_{m_0 +2n_0 +1} .
\end{aligned}
\end{gather*}
Plugging these equalities into $f_{2(m_0+n_0)}=f_0$, we obtain 
\begin{equation}
f_0 = e^{2i\kappa (m_0 +n_0 ) } c_\epsilon^+ f_0,
\label{S4_eq_trajectory00}
\end{equation}
This shows that \eqref{QC2}  is a necessary condition to have a solution $u$ with $u_\leftarrow(0,0)=f_0\neq0$. 
Moreover, under the condition \eqref{EVorRes}, $2(m_0+n_0)$-roots of \eqref{QC2} is real, and $u$ defined by \eqref{EFalongCT} is an associated eigenfunction. 

When the condition \eqref{EVorRes} is false, at least one of the four factors of $c_\epsilon^+$ has its modulus less than 1. For example, let us suppose that $\left|c_{\epsilon,\uparrow,\leftarrow}(0,0) \right|<1$. Then $u$ also satisfies 
\begin{equation*}
e^{-i\kappa}u_\leftarrow(-1,0)=c_{\epsilon,\leftarrow,\leftarrow}(0,0) u_{\leftarrow}(0,0),\quad
e^{-i\kappa}u_\downarrow(0,-1)=c_{\epsilon,\downarrow,\leftarrow}(0,0)u_{\leftarrow}(0,0).
\end{equation*}
Note that under the condition \eqref{S4_eq_Ueassumption3}, $u_\rightarrow(1,0)$ is independent of $u_{\leftarrow}(0,0)$. It follows inductively that
\begin{equation*}
u_\leftarrow(-N,0)=e^{iN\kappa}c_{\epsilon,\leftarrow,\leftarrow}(0,0)f_0,\quad
u_\downarrow(0,-N)=e^{iN\kappa}c_{\epsilon,\downarrow,\leftarrow}(0,0)f_0,
\end{equation*}
for $N=1,2,\ldots$. By a symmetric argument for other corners, we obtain a resonant state. \qed

\medskip

The above proposition with a symmetric argument along $\Phi_-$ shows the following instability of eigenvalues under small perturbations. 

\begin{cor}
For any $\epsilon\in(0,1]$, and any $N\in \{0,2(n_0+m_0),4(n_0+m_0)\}$, one can construct a QW $U_\epsilon$ having $N$ eigenvalues. However, the sum of the numbers of eigenvalues and resonances is always $4(n_0+m_0)$.
\label{cor-Mainof4}
\end{cor}

Proof.
It suffices to show that there is no eigenfunctions or resonant states other than those we have constructed above. If an outgoing solution $u$ to $U_\epsilon u=e^{-i\kappa}u$ for a $\kappa\in{\bf T}_{\bf C}$ does not vanishes at least one point of the image of $\Phi_\pm$, the other entries of $u$ are automatically determined by the argument used in the proof of Proposition \ref{EVorRes-prop}. Thus, $u$ coincides with one of such eigenfunctions or resonant states constructed above. 
Otherwise, that is, if $u$ vanishes at every point of the image of $\Phi_\pm$, then the same argument as in the proof of Lemma \ref{S3_lem_NTresonance} shows that $u$ is identically vanishing.  \qed

\subsection{Asymptotic distribution of eigenvalues and resonances}
In the previous subsection, we saw that $U_\epsilon$ can be eigenvalue-free however small the perturbation is (Corollary \ref{cor-Mainof4}). Contrary, we show in the present subsection that the eigenvalues of the conjugated operator $U(\theta)$ are stable under small perturbations. The following proposition is the main result of this section. Recall that the eigenvalues of $U_{el}$ are given by \eqref{EV-Uel}.

\begin{prop}
Suppose that $U_{\epsilon } $ satisfies \eqref{S4_eq_Ueassumption1}, \eqref{S4_eq_Ueassumption2}, \eqref{S4_eq_Ueassumption3}. 
Let $\theta\in\mathcal{O}_0^-$. 
Then there exists $\epsilon_0>0$ such that for any $\epsilon\in(0,\epsilon_0]$, $U_\epsilon(\theta)=T(\theta)U_\epsilon T(\theta)^{-1}$ has $4(m_0+n_0)$-eigenvalues as well as $U_{el}$. 
Moreover, there exists $r>0$ such that for each $k=0,1,\ldots, 2(m_0 +n_0 )-1$, there exist $\kappa_{k,\epsilon}^+,$ $\kappa_{k,\epsilon}^-\in{\bf T}_{\bf C}$ such that $e^{-i\kappa_{k,\epsilon}^\pm}\in\sigma_p(U_\epsilon(\theta))$ and  
\begin{equation*}
| \kappa_{k,\epsilon}^\pm - \pi k/(m_0 +n_0 ) | <r\epsilon\quad\text{in }\ {\bf T}_{\bf C}
\end{equation*}
for any $\epsilon\in(0,\epsilon_0]$.
\label{S4_prop_ecperturbation}
\end{prop}


Proof. Under the assumption \eqref{S4_eq_Ueassumption2}, we have
\begin{equation*}
c_\epsilon^\pm=(1+O(\epsilon))^4=1+O(\epsilon).
\end{equation*}
Note that $c_\epsilon^+$ is defined by \eqref{EVorRes}, and $c_\epsilon^-$ is by
\begin{equation*}
c_\epsilon^-=c_{\epsilon,\rightarrow,\downarrow}(0,0) c_{\epsilon, \uparrow,\rightarrow}(m_0,0) c_{\epsilon,\leftarrow,\uparrow}(m_0,n_0) c_{\epsilon,\downarrow,\leftarrow}(0,n_0).
\end{equation*}
Recalling the quantization condition \eqref{QC2}, each $\kappa\in{\bf T}_{\bf C}$ satisfying one of
\begin{equation*}
e^{2i(m_0+n_0)\kappa}=c_\epsilon^+=1+O(\epsilon)\quad\text{or }\ 
e^{2i(m_0+n_0)\kappa}=c_\epsilon^-=1+O(\epsilon)
\end{equation*}
is a resonance or an eigenvalue of $U_\epsilon$. Consequently, for each $\kappa$, there exists $k\in\{0,1,\ldots,2(m_0+n_0)-1\}$ such that
\begin{equation*}
\kappa=\frac{\pi k}{m_0+n_0}+O(\epsilon).
\end{equation*}
This shows in particular that $\kappa\in\mathcal{O}^+_{\mathrm{Im} \, \theta}$ (thus $e^{-i\kappa}\in\sigma_p(U(\theta))$) for sufficiently small $\epsilon>0$.
\qed

\section{Shape resonance model for QW}
\subsection{Shape resonance model}
In this section, we redefine the QW $ U_{\epsilon} = SC_{\epsilon} $ as a perturbation of $ U_{np}=SC_{np} = U_i \oplus U_e$ which has been introduced in Subsection 3.2.
Suppose that $ C_{\epsilon} (x)\in \mathrm{U} (4)$ satisfies the assumption
\begin{gather}
C_{\epsilon} (x)= C_{np} (x) \quad \text{for} \quad x\in {\bf Z}^2 \setminus K , \label{S5_eq_shaperesonancemodel} \\
|C_{\epsilon} (x) -C_{np} (x)| _{\infty} < \epsilon \quad \text{for} \quad x\in K, \label{S5_eq_shaperesonancemodel3} 
\end{gather}
for $\epsilon \in(0,1]$.
As a consequence, the estimate
\begin{equation}
\| U_{\epsilon} - U_{np} \| _{{\bf B} (\mathcal{H})} =O(\epsilon ),
\label{S5_eq_shaperesonancemodel4}
\end{equation}
follows.

\subsection{Resolvent estimate}
We fix $\theta \in \mathcal{O}^-_0$. In order to study the resonances, we consider the difference 
$$
T_{\epsilon} (\kappa , \theta )= R_{\epsilon} (\kappa , \theta)-R_{np} (\kappa , \theta ),
$$
of the resolvent operators 
$$
R_{\epsilon} (\kappa , \theta ) = (U_{\epsilon} (\theta)-e^{-i\kappa} )^{-1} , \quad R_{np} ( \kappa , \theta )=R_i (\kappa , \theta ) \oplus R_e (\kappa , \theta ),
$$
where we put $ U_{\epsilon} (\theta )= T(\theta ) U_{\epsilon} T(\theta )^{-1} $ and $ U_{np} (\theta )= T(\theta ) U_{np} T(\theta )^{-1} $.

\begin{lemma}
For $\kappa \in \mathcal{O}^+ _{\mathrm{Im} \, \theta} $ such that $ R_{\epsilon} ( \kappa , \theta )$ and $R_{np} (\kappa , \theta )$ are well-defined, we have
$$
T_{\epsilon} (\kappa,\theta )= R_{\epsilon} (\kappa , \theta )Q_{\epsilon} (\theta ) R_{np} (\kappa , \theta )= R_{np} (\kappa,\theta )Q_{\epsilon} (\theta ) R_{\epsilon} (\kappa , \theta ),
$$
where $Q_{\epsilon} ( \theta )= U_{np} ( \theta )-U_{\epsilon} (\theta )$.
\label{S5_lem_resolventeq_shaperes}
\end{lemma}

Proof.
We put 
$v=R_{np} (\kappa,\theta)f$.
The first equality of the lemma follows from the computation
\begin{gather*}
\begin{split}
(T_{\epsilon} (\kappa , \theta)f,f)_{\mathcal{H}} &= (f,R_{\epsilon} (\kappa,\theta)^* f )_{\mathcal{H}} -(v,f)_{\mathcal{H}} \\ 
&= ((U _{np} (\theta)-e^{-i\kappa})R_{np} (\kappa,\theta)f,R_{\epsilon} (\kappa,\theta)^* f)_{\mathcal{H}} \\ 
&\quad - (v,(U_{\epsilon} (\theta) ^* -\overline{e^{-i\kappa}} ) R_{\epsilon} (\kappa , \theta )^* f )_{\mathcal{H}} ,
\end{split}
\end{gather*}
for any $f\in \mathcal{H}$.
The second equality can be proven in a symmetric argument.
\qed

\medskip

Now we take an eigenvalue $e^{-i\mu_0} \in \sigma_p (U_{np} )$.
Suppose that $\kappa \in {\bf T} _{{\bf C}} $ varies on the counterclockwise loop $\mathcal{L} _{\epsilon,s} ( \mu_0 )$ which has been introduced in Subsection 3.2.
Taking sufficiently small $\epsilon >0$, we assume that there is no other $\mu \in {\bf T} $ such that $e^{-i\mu} \in \sigma_p (U_{np} )$ inside the loop $\mathcal{L} _{\epsilon,s} ( \mu_0 )$.

\begin{lemma}
Fix $s\in(0,1/2]$. 
There exists $\epsilon_0>0$ such that  
$$
\| T_{\epsilon} (\kappa,\theta) \| _{{\bf B} (\mathcal{H})}
$$ 
is uniformly bounded for $\kappa \in \mathcal{L}_{\epsilon,s} (\mu_0 )$ and for $\epsilon\in(0,\epsilon_0]$. 
\label{S5_lem_estimateT}
\end{lemma}

Proof.
Lemma \ref{S3_lem_complexUi} implies
\begin{equation}
\| R_i (\kappa,\theta) \| _{{\bf B} (\mathcal{H}_i)} =O(\epsilon^{-s}), 
\label{S5_eq_estimateT00}
\end{equation}
as well as Lemma \ref{S3_lem_estimateRe} shows the existence of a constant $\gamma >0$ such that 
\begin{equation}
\| R_e (\kappa,\theta) \| _{{\bf B} (\mathcal{H}_e)} \leq \gamma .
\label{S5_eq_estimateT01}
\end{equation}
Recall the operator of restriction $J_e : \mathcal{H} \to \mathcal{H}_e$ which has been used in the proof of Lemma \ref{S3_lem_estimateRe}.
We define $J_i : \mathcal{H}\to \mathcal{H}_i$ in the similar way.
We give the following estimates :
\begin{gather}
\begin{split}
&\| J_i R_{\epsilon} (\kappa,\theta) \| _{{\bf B} (\mathcal{H};\mathcal{H}_i )} \\
&=\| R_i (\kappa , \theta )(U_i (\theta)-e^{-i\kappa} ) J_i R_{\epsilon} (\kappa , \theta) \| _{{\bf B} (\mathcal{H};\mathcal{H}_i)} \\
&\leq \| R_i (\kappa , \theta ) \| _{{\bf B} (\mathcal{H}_i )} \| 1+(U_i (\theta)J_i - J_i U_{\epsilon} (\theta) )R_{\epsilon} (\kappa , \theta ) \| _{{\bf B} (\mathcal{H};\mathcal{H}_i )} \\
&\leq \| R_i (\kappa , \theta ) \| _{{\bf B} (\mathcal{H}_i )} \left( 1+ \| J_i (U_i (\theta) -  U_{\epsilon} (\theta) )R_{\epsilon} (\kappa , \theta ) \| _{{\bf B} (\mathcal{H};\mathcal{H}_i )} \right),
\end{split}
\label{S5_eq_estimateT02}
\end{gather}
and similarly
\begin{gather}
\begin{split}
&\| J_e R_{\epsilon} (\kappa,\theta) \| _{{\bf B} (\mathcal{H};\mathcal{H}_e )} \\
&\leq \| R_e (\kappa , \theta ) \| _{{\bf B} (\mathcal{H}_e )} \left( 1+ \| J_e (U_e (\theta) -  U_{\epsilon} (\theta) )R_{\epsilon} (\kappa , \theta ) \| _{{\bf B} (\mathcal{H};\mathcal{H}_e )} \right) .
\end{split}
\label{S5_eq_estimateT03}
\end{gather}
By the triangular inequality with  (\ref{S5_eq_estimateT00})-(\ref{S5_eq_estimateT03}), we obtain
\begin{gather*}
\begin{split}
&\| R_{\epsilon} (\kappa,\theta) \| _{{\bf B} (\mathcal{H})} \leq \| J_i R_{\epsilon} (\kappa,\theta) \| _{{\bf B} (\mathcal{H};\mathcal{H}_i)} + \| J_e R_{\epsilon} (\kappa,\theta) \| _{{\bf B} (\mathcal{H};\mathcal{H}_e)} \\
&\leq O(\epsilon^{-s} ) \left( 1+O(\epsilon) \| R_{\epsilon} (\kappa,\theta) \| _{{\bf B} (\mathcal{H} )}  \right) + \gamma \left(  1+O(\epsilon)\| R_{\epsilon} (\kappa,\theta) \|_{{\bf B} (\mathcal{H})} \right) \\
&= O(\epsilon^{-s}) +O(\epsilon^{1-s}) \| R_{\epsilon} (\kappa,\theta) \| _{{\bf B} (\mathcal{H} )} .
\end{split}
\end{gather*}
Since $\epsilon>0$ is sufficiently small, we can show
\begin{equation}
\| R_{\epsilon} (\kappa , \theta ) \| _{{\bf B} (\mathcal{H} )} \leq (1+O(\epsilon^{1-s}))^{-1} O(\epsilon^{-s} ) \leq O(\epsilon^{-s}) ,
\label{S5_eq_estimateT04}
\end{equation}
for $s\in (0,1)$.
Now we apply Lemma \ref{S5_lem_resolventeq_shaperes} in order to show
\begin{gather*}
\begin{split}
\| T_{\epsilon} (\kappa,\theta) \| _{{\bf B} (\mathcal{H} )} & \leq \| R_{\epsilon} (\kappa,\theta) \| _{{\bf B} (\mathcal{H})} \| Q_{\epsilon} (\theta) \| _{{\bf B} (\mathcal{H})} \| R_{np} (\kappa,\theta)\| _{{\bf B} (\mathcal{H})} \\
&=  O(\epsilon^{-s})O(\epsilon)O(\epsilon^{-s})= O(\epsilon^{1-2s})=O(1),
\end{split}
\end{gather*}
for $s\in (0,1/2] $.
\qed


\subsection{Existence of resonance}
The eigenvalues of $U_{np}$ are unstable under the perturbation like those of $U_{el}$ where we discussed in the previous section. 
For example, in view of a sufficient condition (C), if neither $\det [c_{np,j,k} (x)] _{j,k\in \{ \leftarrow , \downarrow \} }$ nor $\det [c_{np,j,k} (x)] _{j,k\in \{ \rightarrow , \uparrow \} } $ vanishes for any $x\in \Omega^i\setminus K$, we can take $U_\epsilon$ such that it is eigenvalue-free for any small $\epsilon\in(0,1)$. 
Here, $c_{np,j,k}(x)$ is the $(j,k)$-entry of the coin matrix $C_{np}(x)$. 
As a concrete example, we derive the following case.
In view of the definition of $ C_{np} (x)$, it is a unitary matrix of the form
$$
C_{np} (x)= \begin{bmatrix} 0 & 1 & 0 & 0 \\
 c_{\rightarrow,\leftarrow } (x) & 0 & c_{\rightarrow,\downarrow} (x) & c_{\rightarrow,\uparrow} (x) \\
 c_{\downarrow,\leftarrow} (x) & 0 & c_{\downarrow,\downarrow} (x) & c_{\downarrow,\uparrow} (x) \\
 c_{\uparrow,\leftarrow} (x) & 0 & c_{\uparrow,\downarrow} (x) & c_{\uparrow,\uparrow} (x)   \end{bmatrix} , \quad x\in K_1^- ,
$$
$$
C_{np} (x)= \begin{bmatrix} 0 & c_{\leftarrow,\rightarrow} (x) & c_{\leftarrow,\downarrow} (x) & c_{\leftarrow,\uparrow} (x) \\ 
1 & 0 & 0 & 0  \\ 
0 & c_{\downarrow,\rightarrow} (x) & c_{\downarrow,\downarrow} (x) & c_{\downarrow,\uparrow} (x) \\
 0 & c_{\uparrow,\rightarrow} (x) & c_{\uparrow,\downarrow} (x) & c_{\uparrow,\uparrow} (x)   \end{bmatrix} , \quad x\in K_1^+ ,
$$
$$
C_{np} (x)= \begin{bmatrix} c_{\leftarrow,\leftarrow} (x) & c_{\leftarrow,\rightarrow} (x) & c_{\leftarrow,\downarrow} (x) & 0 \\ 
c_{\rightarrow,\leftarrow } (x) & c_{\rightarrow,\rightarrow} (x) & c_{\rightarrow,\downarrow} (x) & 0 \\ 
 0 & 0 & 0 & 1 \\  
 c_{\uparrow,\leftarrow} (x) & c_{\uparrow,\rightarrow} (x) & c_{\uparrow,\downarrow} (x) & 0 \end{bmatrix} , \quad x\in K_2^- ,
$$
$$
C_{np} (x)= \begin{bmatrix} c_{\leftarrow,\leftarrow} (x) & c_{\leftarrow,\rightarrow} (x) & 0 & c_{\leftarrow,\uparrow} (x) \\ 
c_{\rightarrow,\leftarrow } (x) & c_{\rightarrow,\rightarrow} (x) & 0 & c_{\rightarrow,\uparrow} (x) \\    
 c_{\downarrow,\leftarrow} (x) & c_{\downarrow,\rightarrow} (x) & 0 & c_{\downarrow,\uparrow} (x)  \\
0 & 0 & 1 & 0  \end{bmatrix} , \quad x\in K_2^+ ,
$$
where $c_{j,k} (x)= c_{np,j,k} (x)$ for $j,k\in \{ \leftarrow, \rightarrow, \downarrow, \uparrow\} $.
As the matrix $ C_{\epsilon} (x)$ for $x\in K$, we can take 
$$
C_{\epsilon} (x)= \begin{bmatrix} \epsilon & \sqrt{1-\epsilon^2} & 0 & 0 \\
 \sqrt{1-\epsilon^2} c_{\rightarrow,\leftarrow } (x) & -\epsilon c_{\rightarrow,\leftarrow} (x) & c_{\rightarrow,\downarrow} (x) & c_{\rightarrow,\uparrow} (x) \\
 \sqrt{1-\epsilon^2} c_{\downarrow,\leftarrow} (x) & -\epsilon c_{\downarrow,\leftarrow} (x) & c_{\downarrow,\downarrow} (x) & c_{\downarrow,\uparrow} (x) \\
 \sqrt{1-\epsilon^2} c_{\uparrow,\leftarrow} (x)  & -\epsilon c_{\uparrow,\leftarrow} (x)  & c_{\uparrow,\downarrow} (x) & c_{\uparrow,\uparrow} (x)   \end{bmatrix} , \quad x\in K_1^- ,
$$
$$
C_{\epsilon} (x)= \begin{bmatrix} \epsilon c_{\leftarrow,\rightarrow} (x) & \sqrt{1-\epsilon^2} c_{\leftarrow,\rightarrow} (x)  & c_{\leftarrow,\downarrow} (x) & c_{\leftarrow,\uparrow} (x) \\ 
\sqrt{1-\epsilon^2} & -\epsilon & 0 & 0  \\ 
\epsilon c_{\downarrow,\rightarrow} (x) & \sqrt{1-\epsilon^2} c_{\downarrow,\rightarrow} (x) & c_{\downarrow,\downarrow} (x) & c_{\downarrow,\uparrow} (x) \\
 \epsilon c_{\uparrow,\rightarrow} (x)  &  \sqrt{1-\epsilon^2} c_{\uparrow,\rightarrow} (x)  & c_{\uparrow,\downarrow} (x) & c_{\uparrow,\uparrow} (x)   \end{bmatrix} , \quad x\in K_1^+ ,
$$
$$
C_{\epsilon} (x)= \begin{bmatrix} c_{\leftarrow,\leftarrow} (x) & c_{\leftarrow,\rightarrow} (x) & \sqrt{1-\epsilon^2} c_{\leftarrow,\downarrow} (x) & -\epsilon c_{\leftarrow,\downarrow} (x) \\ 
c_{\rightarrow,\leftarrow } (x) & c_{\rightarrow,\rightarrow} (x) & \sqrt{1-\epsilon^2} c_{\rightarrow,\downarrow} (x) & -\epsilon c_{\rightarrow,\downarrow} (x) \\ 
 0 & 0 & \epsilon & \sqrt{1-\epsilon^2} \\  
 c_{\uparrow,\leftarrow} (x) & c_{\uparrow,\rightarrow} (x) & \sqrt{1-\epsilon^2} c_{\uparrow,\downarrow} (x) & -\epsilon  c_{\uparrow,\downarrow} (x) \end{bmatrix} , \quad x\in K_2^- ,
$$
$$
C_{\epsilon} (x)= \begin{bmatrix} c_{\leftarrow,\leftarrow} (x) & c_{\leftarrow,\rightarrow} (x) & \epsilon c_{\leftarrow,\uparrow} (x)  & \sqrt{ 1-\epsilon^2 } c_{\leftarrow,\uparrow} (x) \\ 
c_{\rightarrow,\leftarrow } (x) & c_{\rightarrow,\rightarrow} (x) & \epsilon c_{\rightarrow,\uparrow} (x)  & \sqrt{ 1-\epsilon^2 } c_{\rightarrow,\uparrow} (x) \\    
 c_{\downarrow,\leftarrow} (x) & c_{\downarrow,\rightarrow} (x) & \epsilon c_{\downarrow,\uparrow} (x) & \sqrt{ 1-\epsilon^2 } c_{\downarrow,\uparrow} (x)  \\
0 & 0 & \sqrt{ 1-\epsilon^2 } & -\epsilon  \end{bmatrix} , \quad x\in K_2^+ .
$$
It is easily checked that the matrix $ C_{\epsilon} (x)$ for every $x\in K$ is unitary.
For the matrix $C_{\epsilon} (x)=[ c_{\epsilon,j,k} (x)]_{j,k\in \{ \leftarrow, \rightarrow, \downarrow, \uparrow\} }$, we have 
\begin{gather*}
\det [ c_{\epsilon,j,k} (x)] _{j,k\in \{ \leftarrow,\downarrow\} } = \left\{ 
\begin{split}
\epsilon c_{\downarrow,\downarrow} (x)&, \quad x\in K_1^- , \\ 
\epsilon ( c_{\leftarrow,\rightarrow} (x) c_{\downarrow,\downarrow} (x)-c_{\leftarrow,\downarrow} (x) c_{\downarrow,\rightarrow} (x)) &, \quad x\in K_1^+ , \\
\epsilon c_{\leftarrow,\leftarrow} (x) &, \quad x\in K_2^- , \\ 
\epsilon (  c_{\leftarrow,\leftarrow} (x)c_{\downarrow,\uparrow} (x)-c_{\leftarrow,\uparrow} (x) c_{\downarrow,\leftarrow} (x)) &, \quad x\in K_2^+ ,
\end{split}
\right.
\end{gather*}
and
\begin{gather*}
\det [ c_{\epsilon,j,k} (x)] _{j,k\in \{ \rightarrow,\uparrow\} } = \left\{ 
\begin{split}
\epsilon (c_{\rightarrow,\uparrow } (x) c_{\uparrow,\leftarrow} (x)- c_{\leftarrow,\rightarrow} (x) c_{\uparrow,\uparrow} (x)) &, \quad x\in K_1^- , \\ 
-\epsilon c_{\uparrow,\uparrow} (x) &, \quad x\in K_1^+ , \\
\epsilon ( c_{\rightarrow,\downarrow} (x) c_{\uparrow,\rightarrow} (x)-c_{\rightarrow,\rightarrow} (x)c_{\uparrow,\downarrow} (x)) &, \quad x\in K_2^- , \\ 
-\epsilon c_{\rightarrow,\rightarrow} (x) &, \quad x\in K_2^+ .
\end{split}
\right.
\end{gather*}
Now suppose that neither $\det [ c_{np,j,k} (x)] _{j,k\in \{ \leftarrow,\downarrow\} }$ nor $\det [ c_{np,j,k} (x)] _{j,k\in \{ \rightarrow,\uparrow\} }$ vanishes for all $x\in \Omega^i \setminus K$.
If $ C_{np} (x)$ for every $x\in K$ satisfies
\begin{gather*}
c_{\downarrow,\downarrow} (x)\not=0 , \quad c_{\rightarrow,\uparrow } (x) c_{\uparrow,\leftarrow} (x)- c_{\leftarrow,\rightarrow} (x) c_{\uparrow,\uparrow} (x) \not= 0 , \quad x\in K_1^- , \\ 
c_{\uparrow,\uparrow} (x)\not= 0 , \quad c_{\leftarrow,\downarrow} (x) c_{\downarrow,\rightarrow} (x)- c_{\leftarrow,\rightarrow} (x) c_{\downarrow,\downarrow} (x) \not= 0 , \quad x\in K_1^+ , \\ 
c_{\leftarrow,\leftarrow} (x)\not= 0, \quad c_{\rightarrow,\downarrow} (x) c_{\uparrow,\rightarrow} (x)- c_{\rightarrow,\rightarrow} (x)c_{\uparrow,\downarrow} (x) \not= 0 ,\quad  x\in K_2^- , \\ 
c_{\rightarrow,\rightarrow} (x)\not= 0 , \quad c_{\leftarrow,\uparrow} (x) c_{\downarrow,\leftarrow} (x)-  c_{\leftarrow,\leftarrow} (x)c_{\downarrow,\uparrow} (x) \not= 0 , \quad x\in K_2^+ ,
\end{gather*}
it follows that $ C_{\epsilon} (x)$ satisfies the condition (C).
One of more trivial cases is given by 
$$
C_{\epsilon} (x)= \begin{bmatrix}
\epsilon & \sqrt{1-\epsilon^2} & 0 & 0 \\ 
\sqrt{1-\epsilon^2} & -\epsilon & 0 & 0 \\
0 & 0 & \epsilon & \sqrt{1-\epsilon^2} \\
0 & 0 & \sqrt{1-\epsilon^2} & -\epsilon \end{bmatrix} , \quad x\in K,
$$
and $ C_{np} (x)= C_{\epsilon} (x)| _{\epsilon =0} $ for $x\in K$.

Now we shall prove the existence result of resonances for $U_{\epsilon}$ near each eigenvalue of $U_{np}$.
We fix $e^{-i\mu_0} \in \sigma_p (U_{np} )$.
In view of (\ref{S2_eq_algebraicmulti_proj}) and (\ref{S2_eq_loop01}), we define the projection operators
\begin{gather*}
\begin{split}
&P_{U_{\epsilon}} (\mu_0 )= \frac{1}{2\pi } \oint _{\mathcal{L} _{\epsilon,s} (\mu_0)} e^{-i\kappa} R_{\epsilon} (\kappa , \theta)d\kappa , \\
&P_{U_{np}} (\mu_0 )= \frac{1}{2\pi } \oint _{\mathcal{L} _{\epsilon,s} (\mu_0)} e^{-i\kappa} R_{np} (\kappa , \theta)d\kappa ,
\end{split}
\end{gather*}
for sufficiently small $\epsilon >0$ and $s\in (0,1/2] $.
Let
$$
\widetilde{P} _{\epsilon} (\mu_0)= P_{U_{\epsilon}} (\mu_0 ) - P_{U_{np} } (\mu_0 )= \frac{1}{2\pi} \oint _{\mathcal{L} _{\epsilon,s} (\mu_0)} e^{-i\kappa} T_{\epsilon} (\kappa , \theta )d\kappa .
$$

\begin{theorem}
Fix $\theta \in \mathcal{O}^-_0 $.
For sufficiently small $\epsilon >0$ and $s\in (0,1/2]$, we have 
\begin{equation*}
\left\|\widetilde{P} _{\epsilon} (\mu_0)\right\|_{{\bf B}(\mathcal{H})}=O(\epsilon^s).
\end{equation*}
\label{S5_thm_existenceres}
\end{theorem}

Proof.
For any $\theta \in \mathcal{O}^-_0$, we apply Lemma \ref{S5_lem_estimateT} which states the uniform boundedness of $\|T_{\epsilon} (\kappa , \theta)\|_{{\bf B}(\mathcal{H})}$ on the loop $\mathcal{L}_{\epsilon,s}(\mu_0)$. Since the length of the integral path, the loop $\mathcal{L}_{\epsilon,s}(\mu_0)$, is of $O(\epsilon^s)$, we obtain the required estimate.
\qed

\medskip

As a consequence, we obtain the following result immediately.
Let $\mu_0\in{\bf T}$ be such that $e^{-i\mu_0} \in \sigma_p ( U_{np} )$, and put
$$
\mathrm{exp} (-i\mathcal{L} _{\epsilon,s} (\mu_0 )) = \{ e^{-i\kappa} \ ; \ \kappa \in \mathcal{L} _{\epsilon,s} (\mu_0 ) \} .
$$

\begin{cor}
Fix $s\in(0,1/2]$ and $\theta\in\mathcal{O}^-_0$.
There exists $\epsilon_0>0$ such that the number of eigenvalues inside the loop $\mathrm{exp} (-i\mathcal{L} _{\epsilon,s} (\mu_0) )$ of $U_\epsilon(\theta)$ coincides with the multiplicity of the eigenvalue $e^{-i\mu_0}$ for any $\epsilon\in(0,\epsilon_0]$, where we count the eigenvalues the same time as their algebraic multiplicity. 
\label{S5_cor_main}
\end{cor}

\textit{Remark.}
Note that the eigenvalues of $U_\epsilon(\theta)$ is either a resonance or an eigenvalue of $U_\epsilon$. Corollary~\ref{S5_cor_main} is an analogue of Proposition \ref{S4_prop_ecperturbation}. 
In view of Proposition \ref{S4_prop_ecperturbation}, the estimate in Corollary \ref{S5_cor_main} does not show the optimal result for the existence of resonances near an eigenvalue of $U_{np}$.
For some models, we expect that Corollary \ref{S5_cor_main} can be improved for $s=1$.

\end{document}